\documentclass[11pt]{amsart}

\usepackage{amssymb}
\usepackage{stmaryrd}
\usepackage[all]{xy}
\usepackage{epsf}

 \newtheorem{theorem}{Theorem}[section]
 \newtheorem{proposition}[theorem]{Proposition}
 
 \newtheorem{corollary}[theorem]{Corollary}
 \newtheorem{conjecture}{Conjecture}
 
 \newtheorem{lemma/definition}[theorem]{Lemma/Definition}
 \newtheorem{proposition/definition}[theorem]{Proposition/Definition}
 \newtheorem{lem/not}[theorem]{Lemma/Notations}

 \newtheorem{example}[theorem]{Example}

 \newtheorem{num}[theorem]{}
 \newcommand{\rmap}{\longrightarrow}

 \newcommand{\comment}[1]{}
\begin{document}

 \title[Birkhoff interpolation with rectangular sets of nodes]
{Uniform Birkhoff interpolation with rectangular sets of nodes }

 % author one information
 \author{Marius Crainic}
\thanks{supported by the Institute for Basic Research in Science (USA), and the Dutch Royal Academy (The Netherlands)}
\author{Nicolae Crainic}
 % Use this \subjclass if you are using amsproc version 2.0 (December 1999).
 % \subjclass[2000]{}

\date{}

% \newpage

 %%% ----------------------------------------------------------------------
 \begin{abstract} 
In this paper we initiate the study of Birkhoff interpolation schemes with main emphasis
on the shape of the set of nodes. We concentrate here on the simplest shapes (``rectangular''). 
The ultimate goal is to obtain a geometrical understanding of the
solvability. We present several regularity criteria, a few conjectures, 
with geometric interpretations, and we illustrate with many examples.
 \end{abstract}
 \maketitle

\tableofcontents

%%%%%%%%%%%%%%%%%%%%%%%%%%%%%%%%%%%%
%%%%%%%%%%%%%%%%%%%%%%%%%%%%%%%%%%%%
%%%%%%%%%%%%%%%%%%%%%%%%%%%%%%%%%%%%
\section{Introduction}             %
\label{Introduction}     %
%%%%%%%%%%%%%%%%%%%%%%%%%%%%%%%%%%%%
%%%%%%%%%%%%%%%%%%%%%%%%%%%%%%%%%%%%
%%%%%%%%%%%%%%%%%%%%%%%%%%%%%%%%%%%%

Apart from its historical importance (see e.g. Newton's interpolation formula,
Gauss' quadrature formula),
interpolation by polynomials still plays a central role in the local
construction of various approximation 
schemes (splines, finite elements, cubature formulas, etc). 
The homogeneous Birkhoff-Hermite problem is one of the most general multivariate polynomial interpolation
problems. In general, it depends on a set $Z\subset\mathbb{R}^n$ of nodes, a ``lower set'' $S$
defining the ``interpolation space'' $\mathcal{P}_{S}$ of polynomials 
(to which the solution is required to belong to), and a set 
$A\subset \mathbb{Z}^{n}_{+}$ of derivatives which appear in 
the interpolation equations. The univariate case ($n=1$) is quite well understood, and it behaves fundamentally different from
the multivariate case \cite{Lo}. For notational simplicity, we restrict here the the bivariate case.

In the literature one finds two types of results. On one hand, there are criteria whose conclusion
holds in the generic case (when the points of $Z$ are in general position). A very good example is the 
work of R.A. Lorentz \cite{Lo}. On the other hand, there are more constructive results which use ideas
from the univariate case to produce explicit solutions for some special classes of interpolation schemes.
A good example is the work of
Gasca and Maeztu \cite{GM} (for more complete list of references, and historical comments, we refer to
the same \cite{GM, Lo}). In such results, although the shape of $Z$ is very special, it is very rare that it is the starting point. More
precisely, the shape comes on the second place: it is the one for which the method works.

Clearly, the patterns of the general problem are far too complex to be understood in an unified manner.
One of the few common features of all approaches is that the shape of $Z$ (and sometimes that of $A$ too) plays
an essential role. On the other hand, when it comes to applications, it often happens that the set of nodes $Z$ is given, 
has a very particular shape (think e.g. of cubature formulas), and there are only special cases/aspects of the problem that 
need to be solved/understood.

In this paper we initiate the study of multivariate Birkhoff interpolation schemes with main emphasis on the shape
of $Z$. After a short discussion on ``shapes'', we will concentrate on the simplest
shapes, namely the rectangular ones (see \ref{rect-shapes}):
\[ Z= \{ (x_i, y_j): 0\leq i\leq p, 0\leq j\leq q \} \]
($p$ and $q$ are non-negative integers).
 Apart from the main definitions, and a short discussion
on ``cartesian shapes'' and their relation with the uniqueness of Birkhoff-Lagrange schemes (see \ref{cartesian-nodes}),
section \ref{Main definitions and constructions} also brings together the main constructions that are relevant to the study of 
rectangular sets of nodes (most importantly, the notion of blow-up presented in \ref{blow}). At the end of the section
we briefly compare the generic case with the rectangular case. In section \ref{Regularity criterias} we present several
regularity criteria which are very useful in examples, and we derive several consequences. Inspired by all the results
and all the examples, we conjecture that lower sets $S$ that are part of a regular scheme with rectangular sets of nodes
must be the result of a blow-up, and we relate this to a stronger version (for rectangular sets of nodes) of the Polya inequalities.
The results presented in this section can also be viewed as a confirmation of the conjectures in many unrelated cases
(depending on the size and shape of $A$).
At the end, we prove that the conjectures hold true in one more case:
$p= q= 1$ (and, this time, for all $A$).  On the other hand, by looking at the pictures associated to regular schemes, one observes another stricking property that regularity seems to imply (a geometric property this time!). Initially, it  was tempting to formulate this property into yet another conjecture, but, to our surpprise, it turned out to be equivalent to the the Polya-type inequalities
mentioned above. This is the content of Theorem \ref{ref-conj2} in Section \ref{Geometric aspects}. In the same section we point out yet another geometric property that regular schemes seem to share
(Conjecture \ref{conj4}, which we use as a guide for constructing interesting examples, cf. e.g. Example \ref{two-of-them}).
And, finally, Section \ref{Examples} contains a large list of examples.

Finally, we would like to point out that, although the results presented here are quite satisfactory for handeling large classes of examples, we feel that the main achievement of this paper is the understanding we gain (e.g. that blow-ups are relevant, that $A$ is geometrically related to $S$, etc). In this direction, we also believe that finding relations with other fields (e.g. algebraic geometry or algebraic topology) would give a new understanding of the conjectures above 
and of the interpolation problem (e.g., note that ``the degrees of multivariate polynomials'' together form an operad, and the blow-up is one of the simplest operations associated to it; see \ref{blow}). 
The existence of a deeper such relation is (yet) another conjecture we would like to adress here.

This paper is part of the second author's PhD dissertation, and it has circulated as a Utrecht University preprint. The present version
follows ``qualified suggestions'' and is considerably shorter (in particular, it skips some of the technical but completely elementary proofs
which will be presented elsewhere \cite{C-Lagrange, C-Polya, C-rectangular, C-lower, C-reduction}). 

 %%%%%%%%%%%%%%%%%%%%%%%%%%%%%%%%%%%%
 %%%%%%%%%%%%%%%%%%%%%%%%%%%%%%%%%%%%
 %%%%%%%%%%%%%%%%%%%%%%%%%%%%%%%%%%%%
 \section{Main definitions and constructions}             %
 \label{Main definitions and constructions}     %
 %%%%%%%%%%%%%%%%%%%%%%%%%%%%%%%%%%%%
 %%%%%%%%%%%%%%%%%%%%%%%%%%%%%%%%%%%%
 %%%%%%%%%%%%%%%%%%%%%%%%%%%%%%%%%%%%

 \begin{num}{\bf Main definitions/terminology:}\rm \ 
In this paper we will use the following notations/terminology:

\begin{itemize}
\item $Z$ denotes a finite set of points in $\mathbb{R}^2$, which plays the role of
 the {\it set of nodes} of the interpolation problem.
\item $A$ denotes a finite set of pairs of non-negative integers, describing the
 (order of the) {\it derivatives} 
 appearing in the interpolation problem. Typical examples are the rectangles
 $R(u, v)$ and the triangles $T(n)$:
 \[ R(u, v)= \{ (i, j)\in \mathbb{Z}_{+}^{2}: i\leq u, j\leq v \},\]
 \[ T(n)= \{ (i, j)\in \mathbb{Z}_{+}^{2}: i+j\leq n \}.\]
\item $S$ is a lower set, that is, a finite subset $S \subset 
 \mathbb{Z}_{+}^{2}$ with the property that
 \[ (u, v)\in S \Longrightarrow R(u, v)\subset S .\] Corresponding to $S$ is the
space 
\[ \mathcal{P}_{S}= \{ P\in \mathbb{R}[x, y]: P= \sum_{(i, j)\in S} a_{i, j} x^i y^j\} ,\]
and $\mathcal{P}_{S}$ plays the role of the {\it interpolation space} (where the
interpolation is taking place). 
\item such triples $(Z, A, S)$ are called {\it interpolation schemes}, and the
 associated problem is: for given
 constants $\{ c_{i, j}(z): (i, j)\in A, z\in Z\}$, find polynomials $P\in
 \mathcal{P}_{S}$ satisfying the equations
 \begin{equation}
\label{interp-eq}
\frac{\partial^{i+j}P}{\partial x^i\partial y^j} (z)= c_{i, j}(z),\  \forall\
 z\in Z,\ (i, j)\in A .
\end{equation}
 One says that $(Z, A, S)$ is {\it solvable} if the problem has solutions for any
 choice of the constants. One says that $(Z, A, S)$ is {\it regular} if the problem has unique solution. 
 The simplest necessary condition for regularity is normality of $(Z, A, S)$, that is,
 $|S|= |Z| |A|$. All schemes in this paper are assumed to be normal.

\item the interpolation equations are in fact a system of linear
equations, and one denotes by $D(Z, A, S)$ its determinant (well defined because of
the normality condition). This is a polynomial
on the coordinates of the nodes ($2n$ such coordinates, where $n= |Z|$). Hence,
fixing $(A, S)$, if $(Z, A, S)$ is regular for one choice of $Z$, then it is regular
for almost all choices of $Z$ (where ``almost all'' refers to the Lebesgue measure on $\mathbb{R}^{2n}$).
In this case one says that $(A, S)$ is {\it almost regular} with respect to sets of $n$ nodes.
\end{itemize}
\end{num}

\begin{num}\label{represent}{\bf Representing lower sets:}\rm \ 
First of all, for any non-negative integers $s$, $n_0\geq n_1\geq \ldots \geq n_s$, 
\[ S_{y}(n_0, \ldots, n_s)= \{ (i, j)\in \mathbb{Z}^{2}: 0\leq i\leq s, 0\leq j\leq n_{i} \} \]
is clearly a lower set. Conversely, any lower set $S$ can be uniquely written in such a way.
Given any set $B\subset \mathbb{Z}^2$, we introduce the notation
\[ B_{y}[\alpha]= \{ \beta: (\alpha, \beta)\in B\} .\]
Then, given a lower set $S$, the associated $s$ is the maximal number $i$ with the property that $S_{y}[i]\neq \emptyset$, 
while $n_i$ is determined by $S_{y}[i]$:
\[ S_{y}[i]= \{ 0, 1, \ldots, n_i\} \]
A similar discussion is obtained by interchanging the role of $x$ and $y$.
Representing $S$ both as $S_{y}(n_0, \ldots, n_s)$ and as $S_{x}(m_{0}, m_{1}, \ldots , m_{t})$,
the relation between the two is:
\[ t= n_{0},\ m_{j}= max \{ i: n_i\geq j\}. \]

\begin{figure}[h]
\begin{center}
         \setlength\epsfxsize{11cm}
        \leavevmode
        \epsfbox{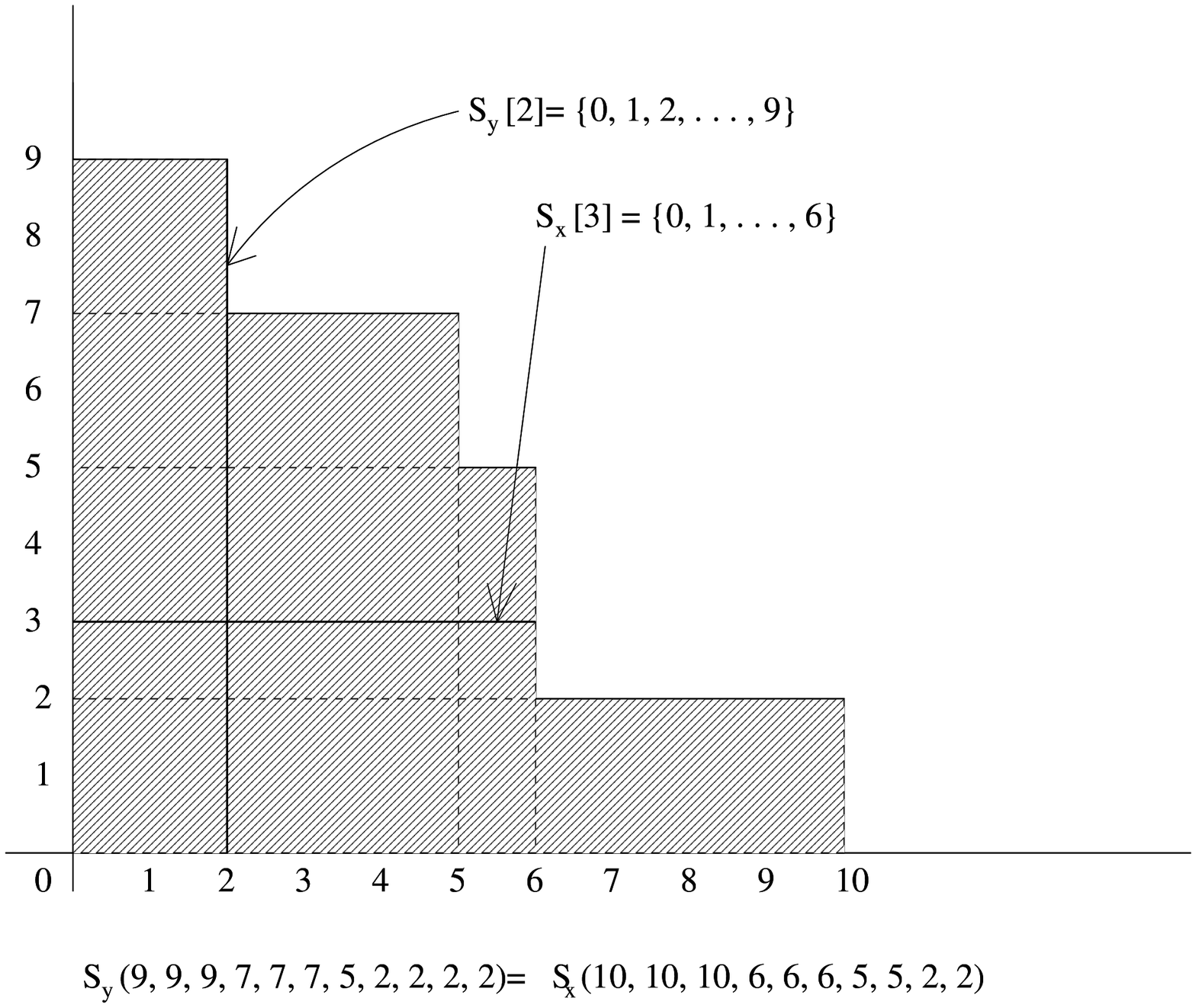}
\end{center}
\caption{\label{exam-fig}}
\end{figure}

Another way of representing lower sets is by making use of its ``exterior boundary points''. Given $S$, 
a point $(u, v)\in S$
is called a {\it boundary point} if $(u+1, v+1)\notin S$. We denote by $\partial S$ the set 
of such points. We consider the following two possibilities:
\begin{enumerate}
\item[(i)] $(u, v+1)\in S$;
\item[(ii)] $(u+1, v)\in S$.
\end{enumerate}

\begin{figure}[h]
\begin{center}
        \epsfxsize=\textwidth
        \leavevmode
        \epsfbox{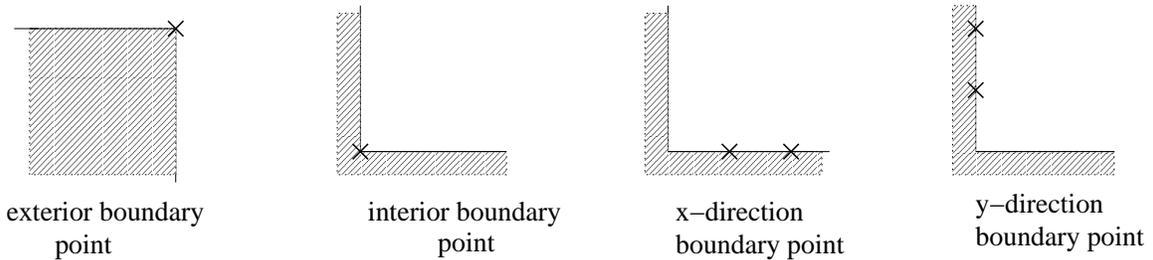}
\end{center}
\caption{Boundary points.}
\label{bound} 
\end{figure}

We denote by $\partial_{e}S$ the set of boundary points $(u, v)$
for which neither of the two conditions above is satisfied (``exterior boundary points''),
by $\partial_{i}S$ the set of those which satisfy both conditions (``interior boundary points''),
by $\partial_{x}S$ the set of those for which only (ii) holds true (``$x$-direction boundary points''), and by
$\partial_{y}S$ the set of those for which only (i) holds true (``$y$-direction boundary points''). 
These four sets form a partition of the boundary $\partial(S)$ of $S$.

The set $\partial_{e}S$ of exterior boundary points determines $S$ uniquely,
since
\[ S= \bigcup_{(u, v)\in  \partial_{e}S} R(u, v) .\]
It is actually useful to label the elements of $\partial_{e}S$ by 
$(a_k, b_1)$, $(a_{k-1}, b_2)$, $\ldots$ , $(a_{1}, b_{k})$, 
indexed so that $a_1<a_2<\ldots < a_k$, and $b_1<b_2<\ldots < b_k$. 

\begin{figure}[h]
\begin{center}
         \setlength\epsfxsize{11cm}
        \leavevmode
        \epsfbox{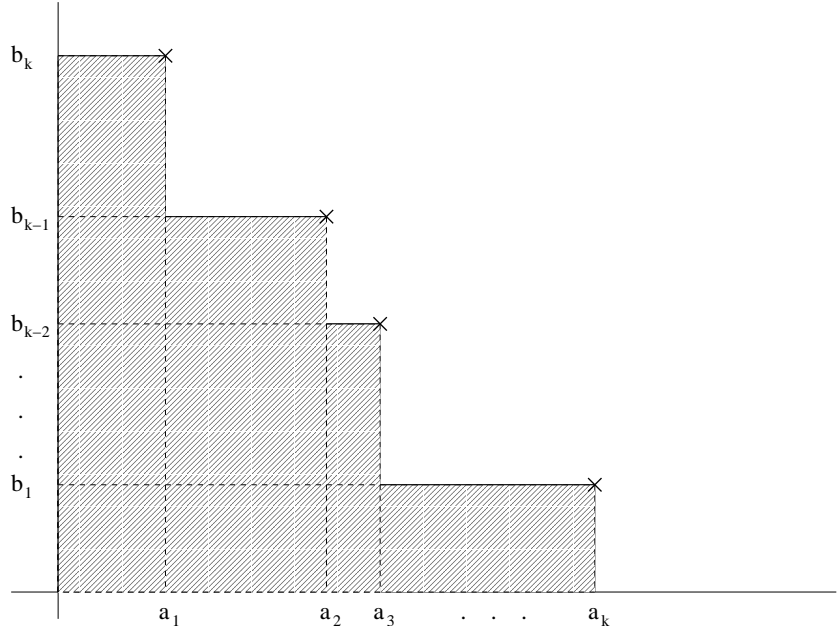}
\end{center}
\caption{\label{edge} Exterior boundary points.}
\end{figure}
Then an element $(u, v)$ is in $S$ if and only if
\[ a_{l-1}< u\leq a_l,\ 0\leq v\leq b_{k- l+ 1}\]
for some $l\in \{ 0, 1, \ldots k\}$. Also,
\[ \partial_{e}S= \{ (a_i, b_j): i+j= k+1\},\  \partial_{i}L= \{ (a_i, b_j): i+j= k\},\]
\[ \partial_{x}S= \{(u, b_{k-j}): 1\leq j\leq k, a_{j}< u< a_{j+1}\}, \]
\[ \partial_{y}S= \{(a_{k-i}, v): 1\leq i\leq k, b_{i}< v< b_{i+1}\}. \]
\end{num}

\begin{num}{\bf Blowing up lowers sets:}
\label{blow}\rm \
Central to Birkhoff interpolation with rectangular sets of nodes is the notion of
blow-up. 
Given $p, q$ non-negative integers, and a lower set $S$, we definer a new lower
set $S^{p, q}$ which is obtained by ``blowing up'' each of the points of $S$ to a
$(p, q)$-rectangle. One may think that a copy of $R(p, q)$ is placed on each of
the points of $S$, and then one pushes these rectangles minimally to the right and upwards
until they become disjoint. For an example with $p= 1$, $q= 2$, see the picture.

\begin{figure}[h]
\label{fig}
\begin{center}
         \setlength\epsfxsize{11cm}
        \leavevmode
        \epsfbox{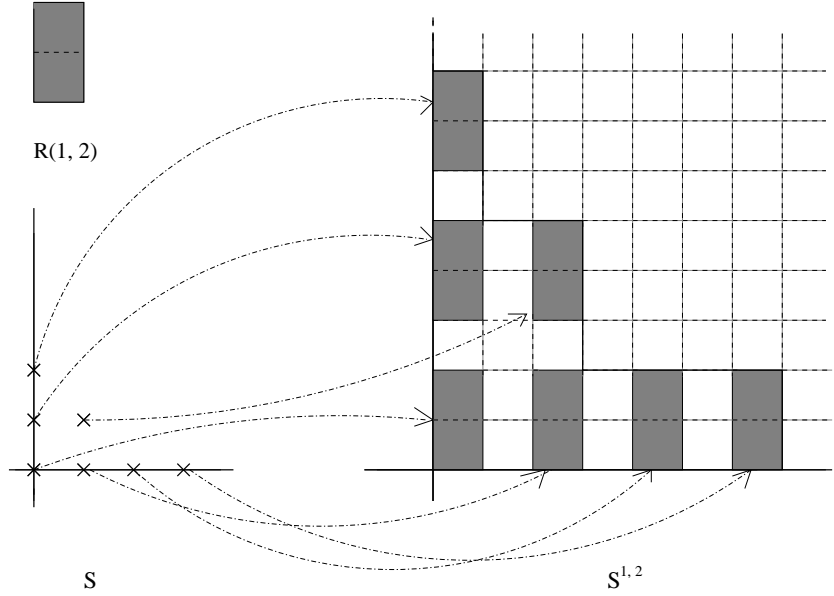}
\end{center}
\caption{\label{blow-fig} Blowing up.}
\end{figure}

More formally, 
\[ S^{p, q}= \{ (\alpha, \beta): ([\frac{\alpha}{p+1}], [\frac{\beta}{q+1}])\in S \}.\]
If $S= S_{y}(n_0, \ldots, n_s)$, then $S^{p, q}= S_{y}(n_{0}^{'}, \ldots, n_{s'}^{'})$, where
\[ k'= (q+1)(k+1)- 1, \ n_{i}^{'}= (p+1)(n_{j}+1)- 1)\ \text{where}\ j= [\frac{i}{q+1}] . \]
In terms of the exterior boundary points, $\partial_{e} S^{p, q}$ consists
of the pairs $((p+1)u, (q+1)v)$ with $(u, v)\in \partial_{e}S$.

It is interesting to point out the structure underlying these operations, structure that becomes
even more important in higher dimensions. Let us denote by $\mathcal{L}(n)$ the set of all
lower sets $L\subset \mathbb{Z}^{n}_{+}$. Then, for any positive integers $n$, $k_1, \ldots , k_n$,
there are ``operations''
\[ \mathcal{L}(n)\times \mathcal{L}(k_1)\times \ldots \times \mathcal{L}(k_n)\rmap \mathcal{L}(k_1+ \ldots + k_n)\]
as follows: given $(S, S_1, \ldots , S_n)$, the result $S'= S\cdot (S_1, \ldots , S_n)$ of this operation is 
the lower set with the property that
\[ \mathcal{P}_{S'}= \{ P(P_1(x_1, \ldots , x_{k_1}), P_2(x_{k_{1}+1}, \ldots, x_{k_{1}+ k_{2}}), \ldots ): P\in \mathcal{P}_{S}, P_i\in \mathcal{P}_{S_i} \} .\]
These operations make $\{ \mathcal{L}(n)\}$ into an operad \cite{Lod}, which we will call {\it the multi-degree operad}.
The reason for the terminology is that lower sets can be viewed as an extension to higher dimension of the ordered set of natural numbers (non-negative integers):
any natural number $n$ defines a lower set $[n]= \{ 0, \ldots , n\}$, and these are all the possible lower sets in dimension one.
In general, an $n$-multivariate polynomials will have an associated lower set playing the role of degree, and $\mathcal{P}_{S}$ will be the space of polynomials ``of degree (at most) $S$''. In particular, lower sets arise naturally also when discussing symbols of multi-deifferantial operators. Note also that the blow-up is one of the simplest operations encoded in the operad $\{ \mathcal{L}(n)\}$:
\[ S^{p, q}= S\cdot ([p], [q]) .\]
Shapes more complicated then the rectangular one would probably require the use of the other operations associated to $\{ \mathcal{L}(n)\}$.
\end{num}

\begin{num}{\bf Particular shapes:}
\label{cartesian-nodes}\rm \ 
We now give a more precise meaning to the ``shape'' of a set of nodes.
First of all, given a lower set $S$, one says that a set of nodes $Z$ is {\it $S$-cartesian}
if it is of type
\[ Z= \{ (x_i, y_j): (i, j)\in S\} ,\]
where the $x_{i}$'s are distinct real numbers, and similarly the $y_{j}$'s. 
We say that $Z$ is cartesian if it is $S$-cartesian for some lower set $S$.

On the other hand, any set of nodes $Z$ induces two lower sets $S_{x}(Z)$, and $S_{y}(Z)$,
which reflect the shape of $Z$. To describe $S_{y}(Z)$, one covers $Z$ 
by lines $l_0, \ldots , l_k$ parallel to the $OY$ axis, and define the numbers
$n_i$ so that on each line $l_i$ there are exactly $n_i+ 1$ points of $Z$. 
We index the lines so that $n_0\geq n_1\geq \ldots \geq n_k$, and we define 
$S_{y}(Z)= S_{y}(n_0, \ldots, n_k)$. The lower set $S_{x}(Z)$ is
defined similarly, by interchanging the role of $x$ and $y$.

The relation between the two is that a set of nodes $Z$ is cartesian if and only if $S_{x}(Z)= S_{y}(Z)$.
In turn, this implies that cartesian sets of nodes show up naturally in relation with the uniqueness of Birkhoff-Lagrange
schemes. Birkhoff-Lagrange schemes are those for which $A$ contains only the origin (and is then suppressed from the notation).
Given any $Z$, there exists at least one lower set $S$ so that the Birkhoff-Lagrange scheme $(Z, S)$ is regular.
However, $S$ is unique if and only if $Z$ is cartesian. For details, see \cite{C-Lagrange}.
\end{num}

\begin{num}\label{rect-shapes}{\bf Rectangular shapes:}\rm \ 
Of special interest for us are the cartesian sets with respect to the rectangles 
$R(p, q)$. Such sets are called {\it $(p, q)$ rectangular} (or just rectangular if we do not want to emphasize 
$p$ and $q$). Hence, $Z$ is {\it $(p, q)$ rectangular} precisely when it is obtained by intersecting
$(p+1)$ vertical lines with $(q+1)$ horizontal lines. Writing
\[ Z= \{ (x_i, y_j): 0\leq i\leq p, 0\leq j\leq q \} \]
where the $x_{0}, \ldots, x_{p}$ are distinct real numbers, and similarly $y_{0}, \ldots, y_{q}$,
the determinant $D(Z, A, S)$ is a polynomial on only $t= p+ q+ 2$
variables (namely the $x_i$'s, and the $y_j$'s). We then say that a pair $(A, S)$ 
is almost regular with respect to $(p, q)$-rectangular sets of nodes if there exists
such a $Z$ so that $(Z, A, S)$ is regular. As in the generic case, this implies
that $(Z, A, S)$ is regular for almost all choices of the $(p, q)$-rectangular set s 
$Z$ (this time, of course, the ``almost'' part refers to the Lebesgue measure on $\mathbb{R}^{t}$).
The regularity is defined similarly. Of course, one can also allow complex nodes, and talk about (almost) 
regularity with respect to $(p, q)$-rectangular sets of complex nodes.
\end{num}

\begin{num}{\bf Examples:}\rm \ 
To get a feeling about the effect that the shape of $Z$ has on regularity, let us point out two cases.
For $A= \{ (0, 0), (1, 1)$, there are no regular schemes whith rectangular sets of nodes, while there are many 
regular schemes with the set of nodes in general position. Consider now the case where 
\[ A= \{ (0, 0), (0, 1) \} .\]
Then one can show that all schemes $(A, S)$ in which $S$ has the property that $|S|= 2n$ and
$S$ contains at most $n$ elements on the $OX$ axis, are almost regular with respect to sets of $n$ nodes.
On the other hand, there is only one scheme which is almost regular with respect to $(p, q)$-rectangular sets 
of nodes, namely $S= R(1, 0)^{p, q}= R(2p+1, q)$. These will become clear 
in the next sections.
\end{num}

\begin{num}\label{Polya}{\bf Polya conditions:} \rm \ 
To compare the multivariate case with the univariate one, and general sets of nodes with rectangular ones, it is interesting
to look at Polya conditions. These are certain algebraic inequalities that are forced by regularity. The classical Polya condition
\cite{Lo} say that, given the regular scheme $(Z, A, S)$ then, for any lower subset $L\subset S$, 
\begin{equation} 
\label{ordpol}
n|L\cap A|\geq |L| ,
\end{equation}
where $n= |Z|$. This follows quite easily by counting the numbers of zero in the determinant associated
to the interpolation equations (\ref{interp-eq}).
Moreover, in the limit case (when equality holds) $(Z, A\cap L, L)$ must be regular too.
This corollary applies to the univariate case as well. Writing
$A= \{ a_0, a_1, \ldots , a_s\}$ with $a_0< a_1< \ldots < a_s$, the Polya conditions become:
\[ a_{i}\leq n\cdot i,\ \ \forall \ 0\leq i\leq s .\]
Moreover, this condition actually insures regularity for almost all sets
of nodes $Z$ \cite{Fe}. More precisely, given $(A, S)$ with $|S|= n|A|$, $(A, S)$
is almost regular if and only if it satisfies the Polya conditions. 
Moreover, if $n=2$, then the Polya conditions are sufficient also for 
regularity. Such properties do not hold in the multivariate case.

When the shape of $Z$ is particular, one can strengthen these inequalities. For instance,
using boundary points (see \ref{represent} above), and by a careful analysis of the zeros in the 
determinant associated to interpolation equations, one can prove \cite{C-Polya} that if $(Z, A, S)$ is regular and $Z$ is $(p, q)$-rectangular, then,
for any lower subset $L\subset S$,
\[ n|A\cap L|\geq |L|+ pq |A\cap \partial L|+ (p+q)|A\cap \partial_{e}L|+
p|A\cap \partial_{y}L|+ q|A\cap \partial_{x}L| .\]
\end{num}

%%%%%%%%%%%%%%%%%%%%%%%%%%%%%%%%%%%%
%%%%%%%%%%%%%%%%%%%%%%%%%%%%%%%%%%%%
%%%%%%%%%%%%%%%%%%%%%%%%%%%%%%%%%%%%
\section{Regularity criteria}             %
\label{Regularity criterias}     %
%%%%%%%%%%%%%%%%%%%%%%%%%%%%%%%%%%%%
%%%%%%%%%%%%%%%%%%%%%%%%%%%%%%%%%%%%
%%%%%%%%%%%%%%%%%%%%%%%%%%%%%%%%%%%%

In this section we bring together several regularity criteria for uniform Birkhoff interpolation schemes
with rectangular sets of nodes, then we present some consequences and then point out a property
that schemes with rectangular sets of nodes seem to satisfy (as suggested by the general criteria, proven
for $|A|\leq 6$, and also for $p= q= 1$). Since the proofs of these criteria, although a bit technical, are completely elementary,
they will be presented elsewhere \cite{C-rectangular, C-lower, C-reduction}. 
Many examples illustrating the usefulness of the criteria will be presented in the last section.

Throughout this section, $(Z, A, S)$ is a uniform Birkhoff scheme,
with a $(p, q)$-rectangular set of nodes
\[ Z= \{ (x_i, y_j): 0\leq i\leq p, 0\leq j\leq q\}.\]
We denote by $A_x$ the elements of $A$ on the $OX$ axis, we also use the similar notations $A_y$, $S_x$, $S_y$,
and we put $Z_x= \{ x_0, \ldots , x_p\}$ (the projection of $Z$ on the $OX$ axis), and similarly $Z_y$.

\begin{proposition}
\label{prop1}
If $(Z, A, S)$ is regular, then
\[ |A|\leq |A_x| |A_y| .\]
Moreover, if the equality holds, then $S$ must be:
\[ S= R(p', q'), \ \ p'= (p+1)|A_x|- 1,\ q'= (q+ 1)|A_y|- 1 .\]
\end{proposition}

This property puts restrictions on the set $A$. For instance, if $mix(A)$ is the
number of mixed derivatives coming from $A$, one immediately gets
\begin{equation}
\label{mized}
\sqrt{mix(A)}\leq \sqrt{|A|}- 1 .
\end{equation}
The next property takes care of the case where one of the lower sets $S_y(A)$ or $S_x(A)$ associated
to $A$ (see \ref{cartesian-nodes}) is rectangular. We then have (for the notation $A_{y}[\alpha]$, see \ref{represent}):

\begin{theorem}
\label{rectangular}
 If $Z$ $A$ has the property 
that $S_y(A)$ is $(s, t)$- rectangular (in particular, if $A$ is $(s, t)$- rectangular), then the scheme
$(Z, A, S)$ is regular if and only if 
\begin{enumerate}
\item[(i)] $S= R(p', q')$, with $p'= (s+1)(p+1)-1, q'= (t+1)(q+1)-1$.
\item[(ii)] The univariate scheme $(Z_{x}, A_{x}, S_{x})$ is regular.
\item[(iii)] All the univariate schemes $(Z_{y}, A_{y}[\alpha], S_{y})$, with $\alpha\in A_x$,
are regular.
\end{enumerate}
\end{theorem}

The next criterion shows the influence that lower subsets of $A$ have on $S$.

\begin{theorem}
\label{lower} 
If $(Z, A, S)$ is solvable (i.e. the interpolation equations have at least one
solution, but not necessarily unique), and if $A$ contains a lower set $R$,
then $S$ must contain the lower set $R^{p, q}$. 
\end{theorem}

In particular, if $A$ is lower, one can completely clarify the situation. 
The outcome makes use of the univariate Hermite polynomials for the set of nodes $Z_x$.
For non-negative integers $u$, $a$, and a node $x_s$, we consider 
\[ H_{a}^{u, s}(x)= \phi_{s}(x)\sum_{k= u}^{a} \frac{(x- x_s)^{k}}{u! (k-u)!} (\frac{1}{\phi_{s}})^{(k-u)} (x_s),\]
where 
\[ \phi_s(x)= (x- x_0)^{a+1} \ldots \widehat{(x- x_s)^{a+1}}\ldots (x- x_p)^{a+1} .\]
We make the convention that $H_{a, s}^{u}(x)= 0$ if $a< u$. Similarly, we denote by $H_{b}^{v, t}$
the univariate Hermite polynomials associated to the set of nodes $Z_y$. 
We also write
\[ \partial_{e} S= \{ (a_k, b_1), \ldots , (a_1, b_k)\},\]
as in \ref{represent}. With these, we have:

\begin{theorem}\label{explicit} If the set $A$ of derivatives is lower, then $(Z, A, S)$ is regular if and only if $S= A^{p, q}$. 

Moreover, in this case, for any $(u, v)\in A$, and any node $(x_s, y_t)\in Z$,
the polynomial
\begin{eqnarray*}
 & H_{u, v}^{s, t}(x, y)  = \\
 & H_{u}^{a_1, s}(x)(H_{v}^{b_k, t}(y)- H_{v}^{b_{k-1}, t}(y))+ H_{u}^{a_2, s}(x)(H_{v}^{b_{k-1}, t}(y)- H_{v}^{b_{k-2}, t}(y))+ \\
 & \ldots + H_{u}^{a_{k-1}, s}(x)(H_{v}^{b_2, t}(y)- H_{v}^{b_1, t}(y))+ H_{u}^{a_k, s}(x)H_{v}^{b_1, t}(y) 
\end{eqnarray*}
is the fundamental interpolation polynomial at the node $(x_i, y_j)$, with respect to the derivative $(u, v)$, i.e.
\begin{equation*}
\frac{\partial^{\alpha+ \beta}H_{u, v}^{s, t}}{\partial x^{\alpha}\partial y^{\beta}}(x_i, y_j)= 
\left \{ \begin{array}{ll}
                     1 & \ \mbox{if}\ \alpha= u, \beta= v, i= s, j= t\\
                     0 & \ \mbox{otherwise}
\end{array}  \right. .
\end{equation*}
\end{theorem}

The next property determines $S$ ``around the coordinate axes''. We use again the notations $S_{y}[\beta]$, $S_{x}[\alpha]$
introduced in \ref{represent}.

\begin{theorem}
\label{axes} 
If $(Z, A, S)$ is a regular UR Birkhoff scheme, where the set of nodes $Z$ is  $(p, q)$-rectangular, then 
\begin{eqnarray}
S_{x}[0]= S_{x}[1]= \ldots= S_{x}[q] & = & \{ 0, 1, \ldots , p'\}, \nonumber \\  
S_{y}[0]= S_{y}[1]= \ldots= S_{y}[p]& = & \{ 0, 1, \ldots , q'\}, \nonumber 
\end{eqnarray}
where $p'= (p+1)|A_x|- 1,\ q'= (q+ 1)|A_y|- 1$. 
\end{theorem}

And, finally, the following shows how one can move or remove points of $A$ on the coordinate axes.
We denote by $S'$ the set obtained from $S$ by removing the lowest 
$(p, q)$-rectangle on its most right
(that is, remove the last $(p+1)$ elements of $S$ from each of the lines 
$y= i$, with $0\leq i\leq q$). Note that $S'$ is lower if and only if $|S_{x}(q+1)|\leq p'-(p+1)$.

\begin{theorem}\label{arrange}\label{replace} 
Let $(Z, A, S)$ be an uniform Birkhoff scheme with rectangular set of nodes. Then:
\begin{enumerate}
\item[(i)] If $(Z, A, S)$ is regular, then the univariate scheme $(Z_x, A_x, S_x)$ must be regular.
\item[(ii)] If $\tilde{A}\subset \mathbb{Z}^2$ is obtained from $A$ by moving some of its elements on the $OX$ axis, then $(Z, \tilde{A}, S)$ is still regular if $(Z_x, \tilde{A}_{x}, S_x)$ is.
\item[(iii)] If $A'$ is obtained from $A$ by removing one of its elements from the $OX$ axis, $S'$ is lower, and  and both univariate schemes $(Z_x, A_{x}, S_{x})$ and $(Z_x, A_{x}', S_{x}')$ are regular, then $(Z, A, S)$ is regular if and only if $(Z, A', S')$ is.
\end{enumerate}
\end{theorem}

There is an obvious variation obtained by interchanging $x$ and $y$. Also, 
there is a particularly good way of changing $A$. Let us denoted by $A_{max}$ 
the set obtained by replacing the elements of $A$ on the $OX$ axis with $\{0, p+1, 2(p+1), \ldots, (p+1)(|A_x|-1)\}$.
This is the ``maximal replacement'' allowed by the univariate Polya condition (see \ref{Polya}). Moreover,
$(A_{max, x}, S_x)$ is automatically regular with respect to all sets of $(p+1)$ nodes
(under the normality condition $|S_x|= (p+1)|A_x|$). In particular, applying first (ii) and then (iii), we deduce that
$(Z, A, S)$ is regular if and only if $(Z_x, A_x, S_x)$ and $(Z, A_{max}', S')$ are.

A large number of examples which show the usefulness of these criteria are presented in Section \ref{Regularity criterias}.
Let us now point out several consequences. First of all, using Theorem \ref{axes} to determine $S$, and
then Theorem \ref{arrange} to reduce the size of the problem (and $A_{max}$ to simplify), we immediately deduce

\begin{corollary}
If $A$ contains no mixed derivative, then $(Z, A, S)$ is regular if and only if the univariate schemes $(Z_x, A_x, S_x)$
and $(Z_y, A_y, S_y)$ are regular and 
\[ S= T^{p, q}, \ \ \text{where}\ \  T= R(|A_x|-1, 0)\cup R(0, |A_y|- 1) .\]
\end{corollary}

The next case (one mixed derivative) is easier to state in terms of almost regularity instead of regularity.

\begin{corollary}
Assume that $A$ contains only one mixed derivative, call it $(\alpha, \beta)$, and let $\{(a_i, 0): 0\leq i\leq s\}$ be the
elements of $A$ on the $OX$ axis, and $\{(0, b_j): 0\leq j\leq t\}$ the ones on the $OY$ axis. Then
$(A, S)$ is almost regular with respect
to $(p, q)$-rectangular sets of nodes if and only if $s, t\geq 1$, 
\[ a_{i}\leq i(p+1), b_{j}\leq j(q+1), \alpha\leq 2p+1, \beta\leq 2q+1\]
(for all $i$ and $j$), and 
\[ S= T^{p, q}, \ \ \text{where}\ \  T= R(|A_x|-1, 0)\cup R(0, |A_y|- 1)\cup R(1, 1) .\]
\end{corollary}

\begin{proof} The conditions on the $a_i$'s and the $b_j$'s are just the univariate Polya conditions (see \ref{Polya}).
On the other hand, if we require that no nontrivial polynomial of type $P= (y- y_0)\ldots (y- y_q)Q(x)\in \mathcal{P}_{S}$
may satisfy the homogeneous interpolation equations of the scheme, we see that $|S_{x}[q+1]|\leq 2p+1$, and, similarly
$|S_{x}[p+1]|\leq 2q+1$. Using also Theorem \ref{axes} and the fact that $S$ is lower, these inequalities
immediately imply the last part of the statement. This allows us to remove elements of $A$ from the coordinate axes 
(use Theorem \ref{arrange} and $A_{max}$) and reduce the problem to $s= 1$ and $t= 1$, $a_1= (p+1)$, $b_1= (q+1)$.
Finally, if the inequalities on $\alpha$ and $\beta$ do not hold, we have $(\alpha, \beta)\ni R(p+1, q+1)$,  and,
using the Polya inequality for $L= R(p+1, q+1)$ (cf. \ref{Polya}) we would obtain
\[ 3n \geq (n+ p+q+3)+ 2pq+ p+ q,\] 
i.e. $3n\geq 3n+1$ which is impossible. Conversely, if $\alpha$ and $\beta$ satisfy these inequalities,
then we can move $(p+1, 0)\in A$ to $(\alpha, 0)$, and then $S_y(A)$ becomes rectangular and we can apply Theorem \ref{rectangular}.
\end{proof}

\begin{corollary}
If $(A, S)$ is almost regular with respect to $(p, q)$-rectangular sets of nodes, and $|A|\leq 6$, 
then $S$ is of type $T^{p, q}$ for some lower set $T$.
\end{corollary}

\begin{proof} Using (\ref{mized}), $A$ will contain at most one mixed derivative except 
for the case where $|A|= 6$ and $A$ contains two mixed derivatives. The first cases
follow from the previous corollary, while the last one follows from the limit case of Proposition \ref{prop1}.
\end{proof}

The results above (see also below, and the section on examples) suggest the following 

\begin{conjecture}
\label{conj1}
If a scheme $(A, S)$ is almost regular with respect to $(p, q)$-rectangular sets of nodes,
then $S= R^{p, q}$ for some lower set $R$.
\end{conjecture}

Let us give a stronger version of this conjecture, which also has a geometric interpretation
that is easier to visualize (see the next section). Given a lower set $S$, we denote by $n_{p, q}(S)$ the
number of elements $(\alpha, \beta)\in S$ with the property that $\alpha$ is divisible
by $(p+1)$, and $\beta$ is divisible by $(q+1)$. The following can be interpreted as
a Polya-type inequality for rectangular sets of nodes.

\begin{conjecture} ($p, q)$-Polya conjecture)
\label{conj2}
If a scheme $(A, S)$ is almost regular with respect to $(p, q)$-rectangular sets of nodes,
then
\begin{equation} 
\label{str-Polya}
|L\cap A|\geq n_{p, q}(L) 
\end{equation}
for all lower sets $L\subset S$.
\end{conjecture}

Let us point out the relation with the Polya condition, as well as the relation between the two conjectures.

\begin{proposition}
Consider the scheme $(A, S)$ with $|S|= n|A|$, where $n= (p+1)(q+1)$.
\begin{enumerate}
\item[(i)] For all lower sets $L$, $n_{p, q}(L)\geq \frac{1}{n}|L|$. In particular,
Conjecture \ref{conj2} is a strengthening of the Polya conditions (\ref{ordpol}).
\item[(ii)] $S= R^{p, q}$ for some lower set $R$ if and only if $n_{p, q}(S)\leq |A|$
(and then equality must hold). In particular, Conjecture \ref{conj2} implies Conjecture \ref{conj1}.
\end{enumerate}
\end{proposition}

\begin{proof}
Given a lower set $L$, we define a new lower set $L_{p, q}$ obtained from $L$ by ``collapsing $(p, q)$-rectangles
to points'':
\[ L_{p, q}= \{ (\alpha, \beta): ((p+1)\alpha, (q+1)\beta)\in L \}.\]
First of all, it is clear that $n_{p, q}(S)= |S_{p, q}|$. Secondly, it is easy
to see that $(L^{p, q})_{p, q}= L$, and $L\subset (L_{p, q})^{p, q}$ for all lower sets 
$L$. Passing to cardinalities in the last inequality, we obtain (i).

On the other hand, one immediately sees that $S_{p, q}$ is the smallest lower set with the property its $(p, q)$-blow up
contains $S$. In particular, if $S= R^{p, q}$ for some lower set $R$, then $R$ must
coincide with $S_{p, q}$. So, the condition that $S= R^{p, q}$ for some lower set $R$ is equivalent to saying
that $S\subset (S_{p, q})^{p, q}$ becomes equality. Passing to cardinalities, 
we have $|A|\leq n_{p, q}(S)$, and the requirement is that equality holds. This proves (ii).
\end{proof}

Regarding (the stronger) Conjecture \ref{conj2}, it is verified in all examples we looked at, and the results of
this section can be used to verify it in many other cases. Let us point out that it is true also
when $p\leq 1$, $q\leq 1$. We state here the case $p= q= 1$, which will be very useful when looking at examples.

\begin{theorem}
\label{case1-1}
Given $A$ and $S$, $p= q= 1$, the following are equivalent:
\begin{enumerate}
\item[(i)] The scheme $(A, S)$ is regular for some $(1, 1)$-rectangular set of nodes $Z$.
\item[(ii)] The scheme $(Z, A, S)$ is regular for any $(1, 1)$-rectangular set of nodes $Z$.
\end{enumerate}
In this case moreover, $(A, S)$ satisfies Conjecture \ref{conj2}
(hence also Conjecture \ref{conj1}).
\end{theorem}

\begin{proof} By changing coordinates
\[ (x, y)\mapsto (2\frac{x- x_0}{x_1-x_0}-1, 2\frac{y- y_0}{y_1-y_0}-1 ,\]
we can always assume 
\[ Z= \{ (-1, -1), (-1, 1), (1, -1), (1, 1)\}, \]
and this proves the equivalence of (i) and (ii). 
Assume now that there exists a lower subset $L\subset S$ such that $n_{1, 1}(L)> |A\cap L|$.
This precisely mens that $R= L_{1, 1}$ (where we use the notations from the previous proof)
satisfies $|R|> |A\cap L|$. 
We now consider the space $\mathcal{P}_{2R}$
of polynomials spanned by $x^{2\alpha}y^{2\beta}$ with $(\alpha, \beta)\in R$.
Note that $\mathcal{P}_{2R}\subset \mathcal{P}_{S}$. 
Next, since the polynomials in $\mathcal{P}_{2R}$ depend on $|R|> |A|$ variables (the coefficients),
we find a non-zero polynomial $P\in \mathcal{P}_{2R}$ such that
\begin{equation}
\label{part} 
\frac{\partial^{i+j}P}{\partial x^{i}\partial y^{j}}(1, 1)= 0, \ \ \forall \ \ (i, j)\in A
\end{equation}
But $P\in \mathcal{P}_{2R}$ implies that $P(-x, -y)= P(-x, y)= P(x, -y)= P(x, y)$. Taking derivatives,
these relations still hold true, at least up to a sign. In particular, (\ref{part}) implies that
$P$ satisfies all the equations (\ref{interp-eq}) corresponding to derivatives coming from $A\cap L$,
and with the constants $c_{i, j}(z)$ equal to zero (the homogeneous equations). On the other hand,
since $L$ is lower and $P\in \mathcal{P}_{L}$, all the derivatives of $P$ coming from $A\setminus L$ are identically
zero. Hence $P\in \mathcal{P}_{2R}\subset \mathcal{P}_{S}$ would be a non-trivial solution of the homogeneous equations
associated to our scheme, which contradicts regularity.
\end{proof}

Note also that, with the same argument as above (using the $(p+1)^{th}$ and the $(q+1)^{th}$ roots of unity)
proves the following

\begin{corollary}
If $(A, S)$ is regular with respect to $(p, q)$-rectangular sets of 
complex nodes, then it satisfies Conjecture \ref{conj2} (hence also of Conjecture \ref{conj1}).
\end{corollary}

%%%%%%%%%%%%%%%%%%%%%%%%%%%%%%%%%%%%
%%%%%%%%%%%%%%%%%%%%%%%%%%%%%%%%%%%%
%%%%%%%%%%%%%%%%%%%%%%%%%%%%%%%%%%%%
\section{Geometric aspects}             %
\label{Geometric aspects}     %
%%%%%%%%%%%%%%%%%%%%%%%%%%%%%%%%%%%%
%%%%%%%%%%%%%%%%%%%%%%%%%%%%%%%%%%%%
%%%%%%%%%%%%%%%%%%%%%%%%%%%%%%%%%%%%

In this section we discuss the geometric interpretations of the inequalities (\ref{str-Polya}), 
which are easy to visualize on the picture (see the next section). 

Roughly speaking, this property says that $S$ is obtainable by attaching a copy of the rectangle $R(p, q)$ to each point 
of $A$, and then moving these rectangles upwards or to the right until they become disjoint. To make this more precise,
we need to use ``shifts'' of $A$ in $S$, which are transformations which move $A$ step by step, at each step one of its elements being moved
upwards or to the right on a new position which is still in $S$ and which is not occupied by any other element
of $A$. We also consider the set 
\[ \mathcal{Z}^{p, q}= \{ ((p+1)i, (q+1)j): i, j-  \text{positive\ integers}\ \} ,\]
which can be viewed as the blow up of the lattice $\mathcal{Z}$ of integral points situated
in the first quadrant.

\begin{theorem}
\label{ref-conj2}
Given $(A, S)$ ($|S|= n|A|$, $n= (p+1)(q+1)$), the following are equivalent
\begin{enumerate}
\item[(i)] $|L\cap A|\geq n_{p, q}(L)$ for all lower sets $L\subset S$.
\item[(ii)] there exists a shift of $A$ in $S$ which moves $A$ to $S\cap \mathcal{Z}^{p, q}$.
\end{enumerate}
\end{theorem}

\begin{proof} 
Before starting the proof, let us fix some notations. Given an element $e\in S$,
we denote by $e^r$ the one step translation of $e$ to the right, and similarly $e^u$ 
(upwards translation). Also, for any lower set $L$, we consider two new lower sets
$\stackrel{\circ}{L}$ and $\overline{L}$ such that 
\[ \stackrel{\circ}{L}\subset L\subset \overline{L},\ n_{p, q}(\stackrel{\circ}{L})= n_{p, q}(L)= n_{p, q}(\overline{L}) ,\]
and such that $\stackrel{\circ}{L}$ is the smallest possible one, while $\overline{L}$ is the largest possible.
Note that $\overline{L}= (L_{p, q})^{p, q}$. It is easy to see that $\overline{L_1\cup L_2}= \overline{L_1}\cup\overline{L_2}$, 
and similarly for intersections and also for the operations $\stackrel{\circ}{L}$, and 
\[ n_{p, q}(L)= \frac{1}{n}|\overline{L}| .\]

We now turn to the proof. 
Assume first that such a shift $\Lambda$ exists, and let $L\subset S$ be a lower set. 
Then the points of $L\cap\mathcal{Z}^{p, q}$ will be obtained by moving some of the points 
of $A$, and those points must come from $A\cap L$. Hence $\Lambda$ defines a bijection between
a subset of $A\cap L$ and $L\cap\mathcal{Z}^{p, q}$, and this proves the desired inequality.
We now keep $p$ and $q$ fixed, we denote by $\mathcal{P}$ the set of
pairs $(A, S)$ that satisfy (i) and have $|S|= n|A|$, 
and we prove that (ii) holds for any $(A, S)\in \mathcal{P}$ by induction on $|S|$. 
The starting point is $|S|= n$. Then $A$ is forced to be $\{(0, 0)\}$ (apply the condition to $L= \{ (0, 0)\}$), 
$S$ is forced to be the rectangle $R(p, q)$ (apply the condition to $L= S$), hence
$A=  S\cap \mathcal{Z}^{p, q}= \{ (0, 0)\}$.
Assume now that $(A, S)$ does satisfy (i),
and the implication has been shown for all pairs $(A', S')$ with $|S'|< |S|$. 
If one of the points of
$A$ can be moved one step upwards, or one step to the right, so that the condition (i)
is still satisfied (by $S$ and the new set $A$), then we perform the move. We repeat this
if still possible. This process
will stop at some point (e.g. because $A$ cannot be moved outside $S$). Hence, we may assume that
$A$ is ``maximal'', in the sense that none of its elements can be moved any further without 
violating (i). We will show that $A= S\cap \mathcal{Z}^{p, q}$. We first prove the following:\\

{\it Claim 1:} If $L\subset S$ is a lower subset with $|A\cap L|= n_{p, q}(L)$,
then either $\overline{L}= S$, or $A\cap L= \mathcal{Z}^{p, q}\cap L$.\\

Proof of the claim: We first show that $(A\cap L, \overline{L})\in \mathcal{P}$. First of all,
$|\overline{L}|= n_{p, q}(L)= |A\cap L|$. Secondly, for any lower set $P\subset \overline{L}$
one has $\overline{P}\subset \overline{L}$, hence $\overline{L\cap P}= \overline{L}\cap \overline{P}= \overline{P}$.
This, and the fact that $(A, S)$ does satisfy (i), imply that
\[ n_{p, q}(P)= n_{p, q}(L\cap P)\leq |A\cap (L\cap P)|= |(A\cap L)\cap P| \]
for all $P\subset \overline{L}$ lower. Hence $(A\cap L, \overline{L})\in \mathcal{P}$. Assume now that 
$\overline{L}\neq S$. Then, by the induction hypothesis, there is a shift that moves $A\cap L$ to
$\mathcal{Z}^{p, q}\cap\overline{L}= \mathcal{Z}^{p, q}\cap L$. Since this shift stays 
inside $L$, it does not touches the elements in $A\setminus L$, hence it can be viewed as a shift $\Lambda$ 
of the entire $A$, which leaves $A\setminus L$ intact. In other words, the image $A'$ of $\Lambda$ is
given by 
\[  A'\cap L= \mathcal{Z}^{p, q}\cap L,\ A'\setminus L= A\setminus L .\]
Now, for any $P\subset S$ lower, one has
\[  \begin{split}
|A'\cap P| & = |A'\cap (P\cap L)|+ |A'\cap (P\setminus L)|\\
           & = |A'\cap (P\cap L)|+ |A\cap (P\setminus L)|\\
           & = |(A'\cap L)\cap (P\cap L)|+ |A\cap (P\cup L)|- |A\cap L|\\
           & \geq n_{p, q}(P\cap L)+ n_{p, q}(P\cup L)- n_{p, q}(L)\\
           & = n_{p, q}(R)
\end{split}  \]
In the inequality above we used that $(A'\cap L, \overline{L})$ and $(A, S)$ satisfy (i),
and that $|A\cap L|= n_{p, q}(L)$.
This shows that $(A', S)\in \mathcal{P}$. From the maximality of $A$, the shift must be
the identity shift, which precisely means that $A\cap L= \mathcal{Z}^{p, q}\cap L$.
This concludes the proof of the claim.\\

{\it Claim 2:} If $e\in A\setminus \mathcal{Z}(p, q)$ is an element with
the property that $e^r\notin A$, then
\[ e\in \partial_{y}(\stackrel{\circ}{S})\cup \partial_{e}(\stackrel{\circ}{S}) .\]

Proof of the claim: Denote by $A_{e}^{r}$ the new set obtained from $A$ by moving $e$ to $e^r$. By the maximality of $A$, 
$(A_{e}^{r}, S)$ does not satisfy (i), hence we find a lower set $L$ such that
$|A_{e}^{r}\cap L|< n_{p, q}(L)$. But $|A_{e}^{r}\cap L|$ is either $|A\cap L|- 1$ (if $e\in L$, $e^r\notin L$),
or $|A\cap L|$ (in the remaining cases). Since $|A\cap L|\geq n_{p, q}(L)$, we must have
\[ |A\cap L|= n_{p, q}(L), \ e\in L, \ e^r\notin L .\]
We now use the previous claim. Since $e\in (A\cap L)\setminus \mathcal{Z}(p, q)$,
we must have $\overline{L}= S$. In turn, this implies that $\stackrel{\circ}{L}= \stackrel{\circ}{S}$, hence, since $e^r\notin L$, we have 
$e^r\notin \stackrel{\circ}{S}$. On the other hand, applying (i) to the largest lower set which does not
contain $e$, we immediately see that $S$ contains at least one element in $\mathcal{Z}^{p, q}$ larger
then $e$. In other words, $e\in \stackrel{\circ}{S}$. But the elements $e\in \stackrel{\circ}{S}$ with $e^r\notin \stackrel{\circ}{S}$
are exactly those those situated on $\partial_{y}(\stackrel{\circ}{S})\cup \partial_{e}(\stackrel{\circ}{S})$.
This concludes the proof of the claim.\\

Clearly, one can replace $e^r$ by $e^u$, and arrive to a
similar conclusion. In particular, if $e\in A\setminus \mathcal{Z}^{p, q}$,
then either $e^r$ or $e^u$ must be in $A$ (otherwise we must have $e\in \partial_e(\stackrel{\circ}{S})$, 
which is a contradiction because the extremal boundary points of $\stackrel{\circ}{S}$ belong to 
$\mathcal{Z}^{p, q}$).

We are now ready to prove that $A= L\cap\mathcal{Z}^{p, q}$. If this is not so, then we
pick up an extremal element $e$ of $A\setminus \mathcal{Z}^{p, q}$. That means that  
$e\in A\setminus \mathcal{Z}^{p, q}$, and $A\setminus \mathcal{Z}^{p, q}$ contains no other element
whose coordinates are greater or equal to the coordinates of $e$. 
Then at least one of the positions $e^u$ or $e^r$ are not in $A$. Otherwise,
since $e$ is extremal, $e^r$ and $e^u$ would both be in $\mathcal{Z}^{p, q}$, and that is
clearly impossible. Hence, combined with the previous remark, one (and only one) 
of these two elements are in $A$. Hence we may assume that
\[ e^r\notin A,\ \ e^u\in A .\]
From Claim 2, and the extremality of $e$, we deduce that 
\[ e^r\in \partial_{y}(\stackrel{\circ}{S})\cup \partial_{e}(\stackrel{\circ}{S}),\ e^u\in \mathcal{Z}^{p,q}.\]
Let $e'$ be the exterior boundary point of $\stackrel{\circ}{S}$ with the property that 
the segment $l= [e, e']$ is parallel to $OY$ (see Fig. ~\ref{fproof}). We apply condition 
(i) to the lower set $S\setminus Q_{e}$, where $Q_{e}$ is as in the picture.
(i.e. $S\setminus Q_{e}$ is the largest lower subset of $S$ which does not contain $e$). 
We immediately get $|l\cap \mathcal{Z}^{p, q}|\geq |l\cap A|$.
Since $e$ is the only point of $A$ situated on $l$ but not on $\mathcal{Z}^{p, q}$,
we find a point $f= (\alpha, \beta)\in A\cap l\cap\mathcal{Z}^{p, q}$ with the property that 
$\tilde{f}= (\alpha+ p, \beta)$ is in $l\cap\mathcal{Z}^{p, q}$ but it is not in $A$.
We then consider the shift that moves $f$ to $\tilde{f}$, and leaves the rest of $A$ unchanged.
The maximality of $A$ implies again the existence of a lower set $L$ with the
property that $|A\cap L|= n_{p, q}(L)$, $f\in L$, $\tilde{f}\notin L$. We can now use Claim 1
above. Since $\tilde{f}\notin L$ and $\tilde{f}\in \mathcal{Z}^{p, q}$,
one cannot have $\overline{L}= S$. Hence we must have $A\cap L= \mathcal{Z}^{p, q}\cap L$.
On the other hand, since $L$ is lower and $e$ sits below $f$, it follows that $e\in A\cap L= \in L$, hence 
$e\in \mathcal{Z}^{p, q}$. This contradicts the choice of $e$, and concludes the proof of the Theorem.  
\begin{figure}[h]
\begin{center}
        \setlength\epsfxsize{8cm}
        \leavevmode
        \epsfbox{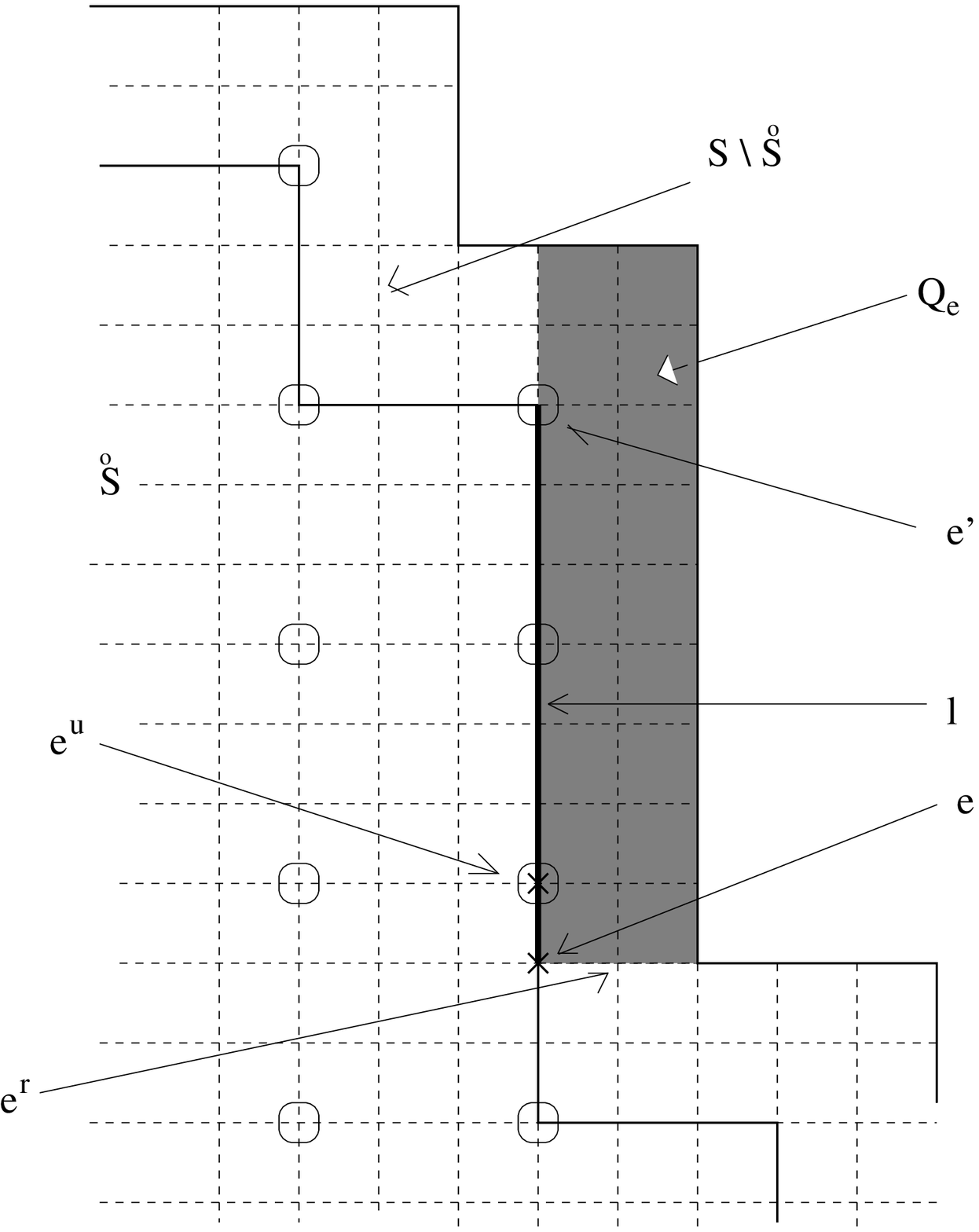}
\end{center}
\caption{}
\label{fproof}
\end{figure}
\end{proof}

Let us point out another geometric property that is suggested by examples. 
This time, we will use {\it inverse shifts}, which are defined exactly like shifts, but moving
downwards and to the left.  

\begin{conjecture}
\label{conj4}
If $(A, S)$ is almost regular with respect to $(p, q)$-rectangular sets of nodes, 
then there is an inverse shift of $A$ in $S$, which moves $A$ to a lower set $R$,
and $S= R^{p, q}$. 
\end{conjecture}

As in Theorem \ref{ref-conj2} (and proven using the same ideas), this geometrical condition is
equivalent to an algebraic condition, which says that $|A\cap L|\leq |R\cap L|$ for all lower sets $L$.
Note also that one should only allow certain type of shifts which do not violate 
certain ``regularity conditions''.  What ``regularity condition'' exactly means
is still to be discovered, but it certainly excludes moving new elements to the axes 
(cf. the results of the previous section). Also, it is tempting to combine the last conjecture and
Theorem \ref{ref-conj2} into a stronger (unifying conjecture) which states the existence of
(certain) shifts of $R$ into $S\cap \mathcal{Z}^{p, q}$, which move $(i, j)$ to $((p+1)i, (q+1)j)$,
and which, at some intermediate step, cover $A$.

Note that in almost all the cases we have considered so far (e.g. no mixed derivatives,
or lower sets of derivatives), moving $A$ backwards to a lower set was possible in only one way.
Hence Conjecture \ref{conj4} would explain the uniqueness of $S$ (proven by us in each case
separately. Also, one can use this conjecture as a guide 
for constructing interesting examples (e.g. where $S$ is not unique, see Example \ref{two-of-them}).

%%%%%%%%%%%%%%%%%%%%%%%%%%%%%%%%%%%%
%%%%%%%%%%%%%%%%%%%%%%%%%%%%%%%%%%%%
%%%%%%%%%%%%%%%%%%%%%%%%%%%%%%%%%%%%
\section{Examples}
\label{Examples}
%%%%%%%%%%%%%%%%%%%%%%%%%%%%%%%%%%%%
%%%%%%%%%%%%%%%%%%%%%%%%%%%%%%%%%%%%
%%%%%%%%%%%%%%%%%%%%%%%%%%%%%%%%%%%%

In this section we present several examples that illustrate the results of the previous sections.
For simplicity, we restrict most of the examples to the case $p= q= 1$, i.e. 
the case of $(1, 1)$-rectangular sets of nodes. One of the simplifications comes from the
fact that, in this case, the notions of regular and almost regular coincide (cf. Theorem \ref{case1-1}).
Hence, unless otherwise specified,  
the term ``regular'' in this section stands for``regular with respect
to $(1, 1)$-rectangular sets of nodes''. Passing to the general case of $(p, q)$-rectangular sets
of nodes require some care on almost regularity versus regularity, but many of the arguments
remain the same.

\begin{example}\rm \
Let $A$ be as in Fig. ~\ref{nonreg}. Then there is no lower set $S$ 
which makes $(A, S)$ into a regular scheme. Assume there is one. 
First of all, it must be $R(5, 5)$, as implied by Proposition \ref{prop1}.
In this case however, the inequality (\ref{str-Polya}) (insured by Theorem \ref{case1-1}) is violated
by the lower set $L$ drawn in the picture. Hence, there is no lower set
$S$ which makes $(A, S)$ regular.
\[ s(Z_{1, 1}, A)= 0 .\]
\begin{figure}[h]
\begin{center}
        \setlength\epsfxsize{11cm}
        \leavevmode
        \epsfbox{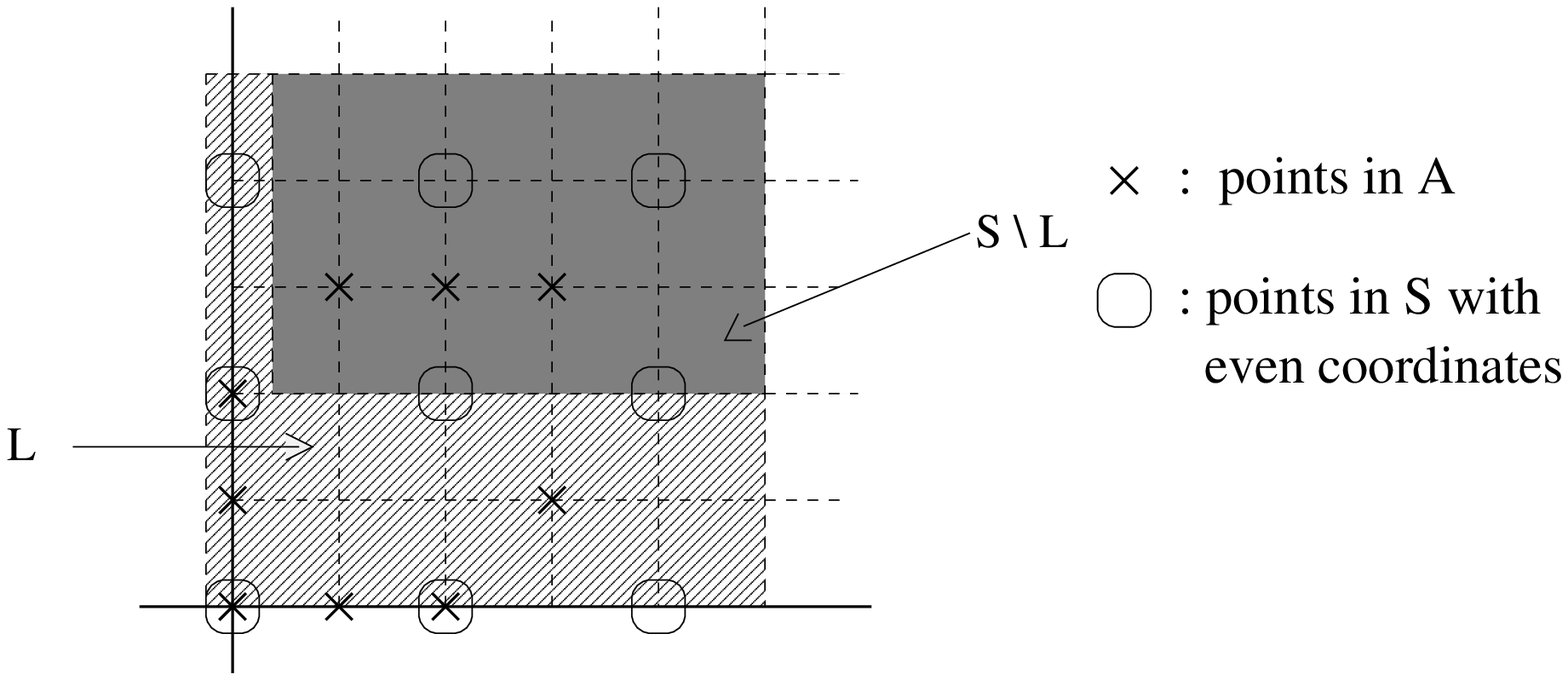}
\end{center}
\caption{}
\label{nonreg}
\end{figure}
\begin{figure}[h]
\begin{center}
        \setlength\epsfxsize{8cm}
        \leavevmode
        \epsfbox{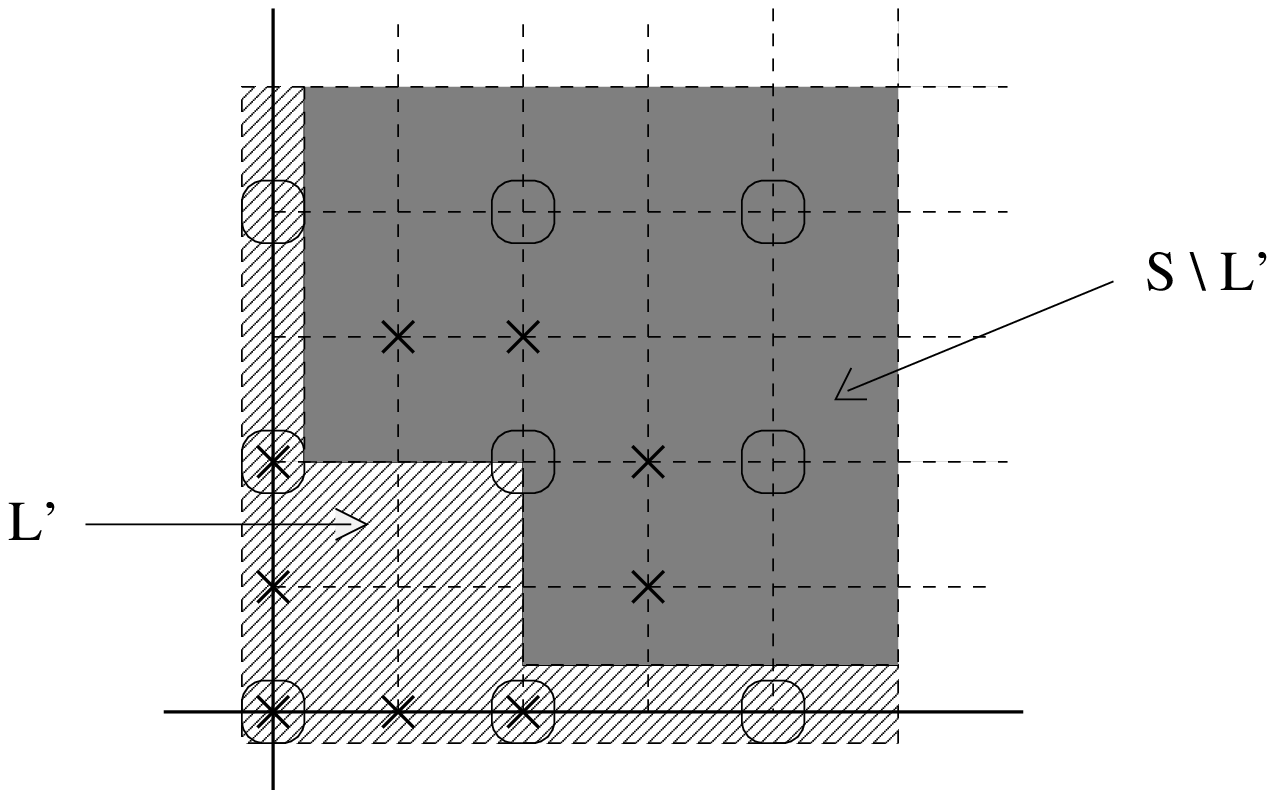}
\end{center}
\caption{}
\label{nonregg}
\end{figure}
Note also that after moving any of the first two points of $A$ situated on the line
$y= 3$ one step downwards, the condition is no longer violated, and
one can actually show that the resulting schemes are regular. 
However, moving the last point on that line (i.e. $(3, 3)$) one step downwards, 
produces a scheme which is still non-regular since the same condition is violated 
(this time by $L'$ shown in Fig. ~\ref{nonregg}).
\end{example}

\begin{example}\rm \
\label{exfig2}
Consider now $A$ as in Fig. ~\ref{reg1}. 
As above, regularity forces $S= R(5, 5)$. This time however, 
the inequality (\ref{str-Polya}) is satisfied. Let us show that
$(A, S)$ is regular. We first remark that $S_{y}(A)$ is $(3, 3)$-rectangular, hence 
we can use Theorem \ref{rectangular} to reduce the regularity of $(A, S)$ to the regularity of
several univariate schemes. In turn, the univariate schemes
are being taking care of by the Polya condition (see \ref{Polya}).
When $p, q\geq 1$ the same arguments apply to conclude that $(A, S)$ 
is almost regular (with respect to $(p, q)$-rectangular sets of nodes) if and only
if $S= R(2p+3, 2q+3)$. 
Moreover, given the set of nodes $Z$, Theorem \ref{rectangular}
rephrases the regularity of $(Z, A, S)$ in terms of the regularity
of certain induced univariate schemes which are easier to handle.
\begin{figure}[h]
\begin{center}
        \setlength\epsfxsize{10cm}
        \leavevmode
        \epsfbox{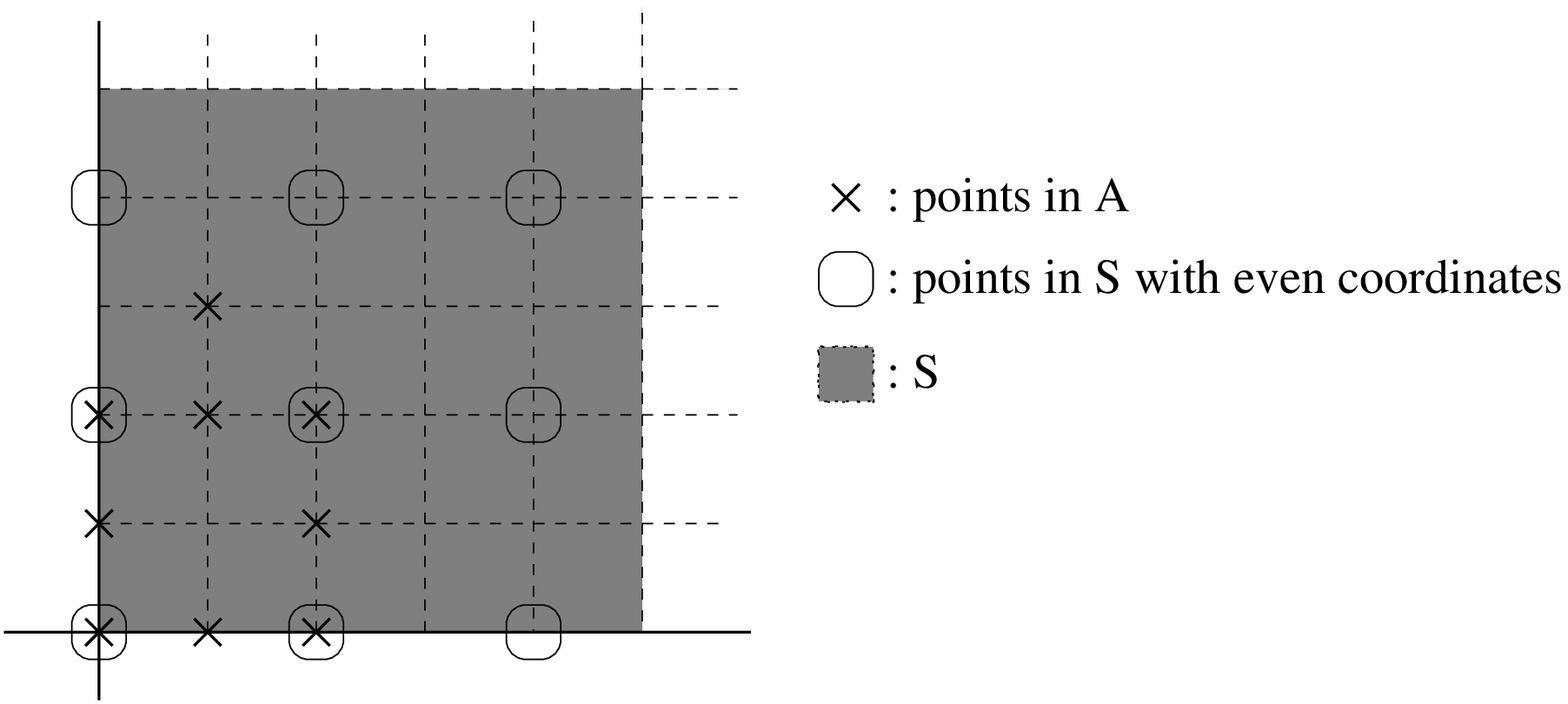}
\end{center}
\caption{}
\label{reg1}
\end{figure} 
\end{example}

\begin{example}\rm \
We consider now the variation of the previous example shown in Fig. ~\ref{reg0}.
One cannot apply Theorem \ref{rectangular} directly, but one can first invoke 
Theorem \ref{replace} to re-arrange the points of $A$ on $OX$
to occupy the first three positions. Then $S_{y}(A)$ becomes rectangular,
and Theorem \ref{replace} can be used. Alternatively,
one can use Theorem \ref{replace} twice (once on each of the axes) to 
and then reduce the problem to the one of the previous example.
\begin{figure}[h]
\begin{center}
        \setlength\epsfxsize{5cm}
        \leavevmode
        \epsfbox{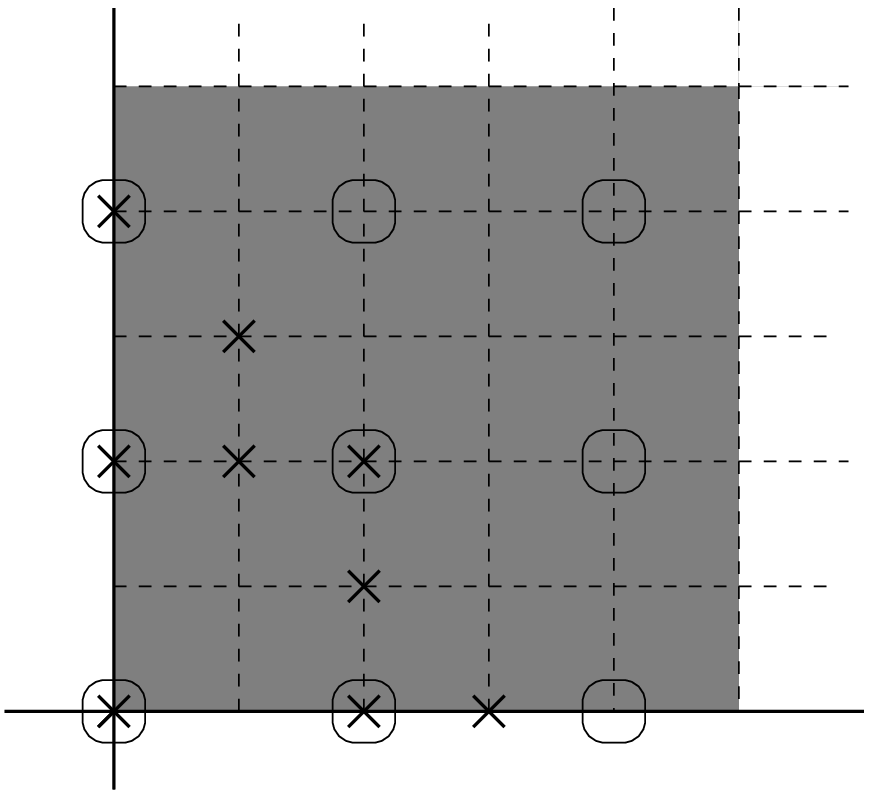}
\end{center}
\caption{}
\label{reg0}
\end{figure} 
\end{example}

\begin{example}\rm \
\label{ten4}
Let $A$ be as in Fig. ~\ref{reg6}.  Using Theorem \ref{axes}, we see that $S$ must contain
$S_0$ shown in the picture, and it must be contained in $R(5, 5)$.
Using Theorem \ref{case1-1} (namely that $S$ is of type $R^{1, 1}$ for some lower set $R$),
we see that $S$ is obtained from $S_0$ together with a copy of the
rectangle $R(1, 1)$. But there are only two ways one can add such a rectangle to $S_0$
to obtain a lower set, and the two possibilities are shown as Case $1$ and Case $2$ in the picture.
In the first case, 
the inequality (\ref{str-Polya}) (insured by Theorem \ref{case1-1})
with $L= S\setminus \{(3, 1)\}$
is not satisfied. The situation is different in the second case, when we obtain a regular scheme. 
To see this, one first uses Theorem \ref{replace}
to remove the last point of $A$ situated on $OY$, and then one treats the remaining scheme
as in Example \ref{exfig2}. Hence, again, there is only one $S$ which makes the scheme regular.
\begin{figure}[h]
\begin{center}
         \setlength\epsfxsize{14cm}
         \leavevmode
         \epsfbox{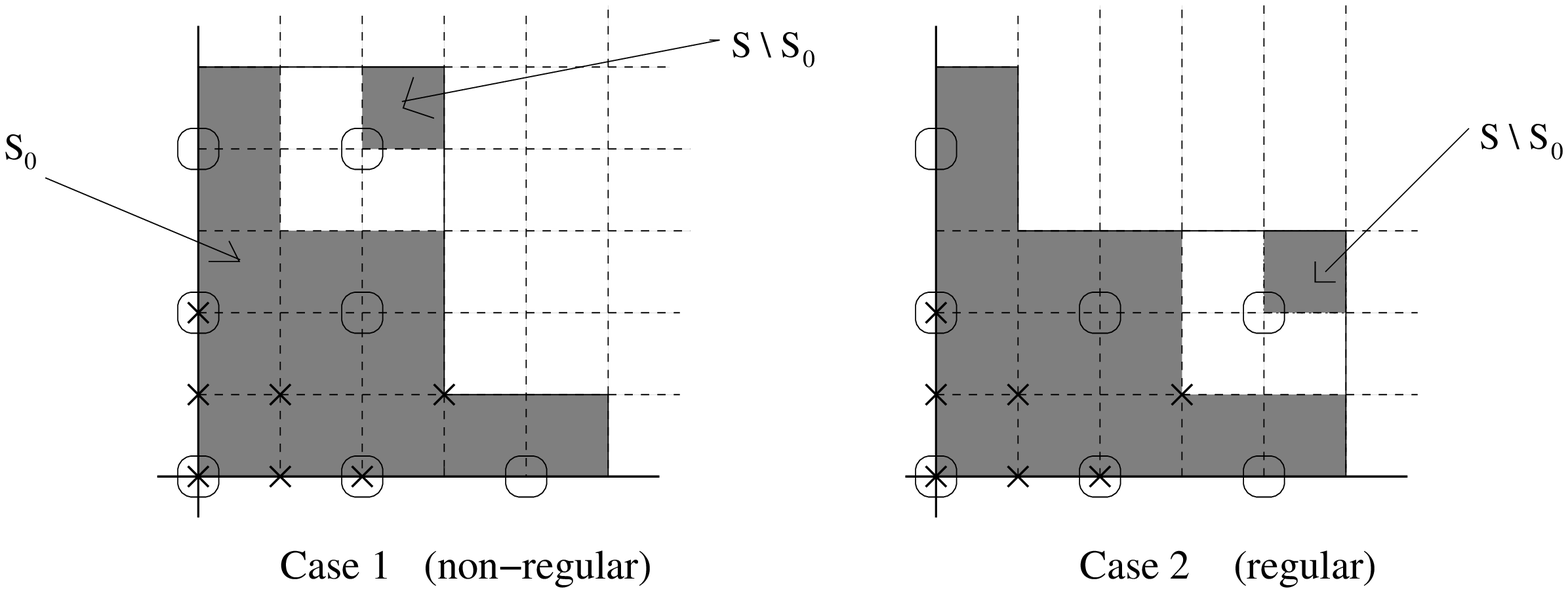}
 \end{center}
 \caption{}
 \label{reg6}
 \end{figure}
\end{example}

\begin{example}\rm \
\label{half}
Let us now describe two examples which show that
the inequalities (\ref{str-Polya}) 
do not imply regularity
with respect to $(1, 1)$-rectangular sets of nodes.
Consider the scheme $(A, S)$ appearing on the left hand side 
of Fig. ~\ref{condb2}. It does satisfy the desired condition,
but it is not regular. To see this, we remark that $A$ can
be obtained from $\tilde{A}$ (see right hand side of
Fig. ~\ref{condb2}) by removing the last element from $OX$.
Hence we can apply Theorem \ref{replace} to $(\tilde{A}, \tilde{S})$
to conclude that it is regular if $(A, S)$ is. But this cannot happen
because the regularity of $(\tilde{A}, \tilde{S})$ and the fact that $\tilde{A}$ 
is lower would imply that $\tilde{S}= \tilde{A}^{1, 1}$ (cf. Theorem \ref{explicit}),
which is not the case.                            
\begin{figure}[h]
\begin{center}
        \setlength\epsfxsize{11cm}
        \leavevmode
        \epsfbox{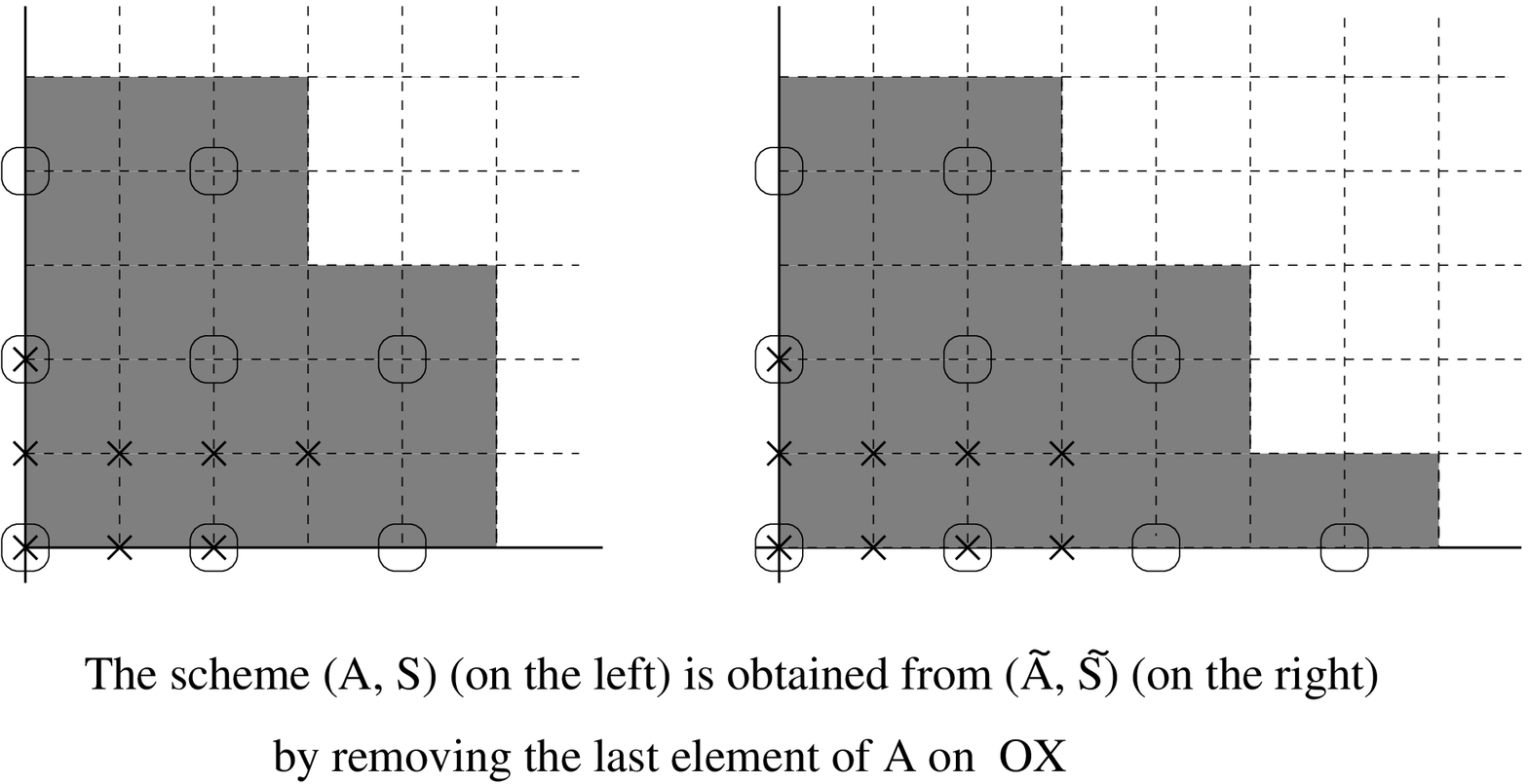}
\end{center}
\caption{}
\label{condb2}
\end{figure} 
We should say here that what causes the non-regularity in this example is another
simple condition that must be satisfied by all regular schemes (and, in this example, it is not): 
the number of points of $A$ on the line $y= 1$ cannot exceed those on the line $y= 0$. 
\end{example}

\begin{example}\rm \
A bit more subtle is the example drawn in Fig. ~\ref{nonreg2}, which still
satisfies the inequalities (\ref{str-Polya}). We advise the reader to try to guess
a ``general regularity condition'' that is broken in this example.
We now give an argument that proves that $(A, S)$ cannot be regular. Assume it is.
We first move the last two elements of $A$ on $OX$ to new positions to get the scheme
$(A', S')$ of Fig. ~\ref{nonreg2}. By Theorem \ref{replace} (and the univariate Polya conditions,
see \ref{Polya}), 
$(A', S')$ is still regular. Now, choose $L$ as in the picture. One has
$|L|= 4 |A\cap L|$, i.e. we are in the limit case of the Polya condition. Hence (see
\ref{Polya}) $(A\cap L, L)$ must be
regular too. But, by the same arguments as
in the example above, $(A\cap L, L)$ cannot be regular.
\begin{figure}[h]
\begin{center}
        \setlength\epsfxsize{11cm}
        \leavevmode
        \epsfbox{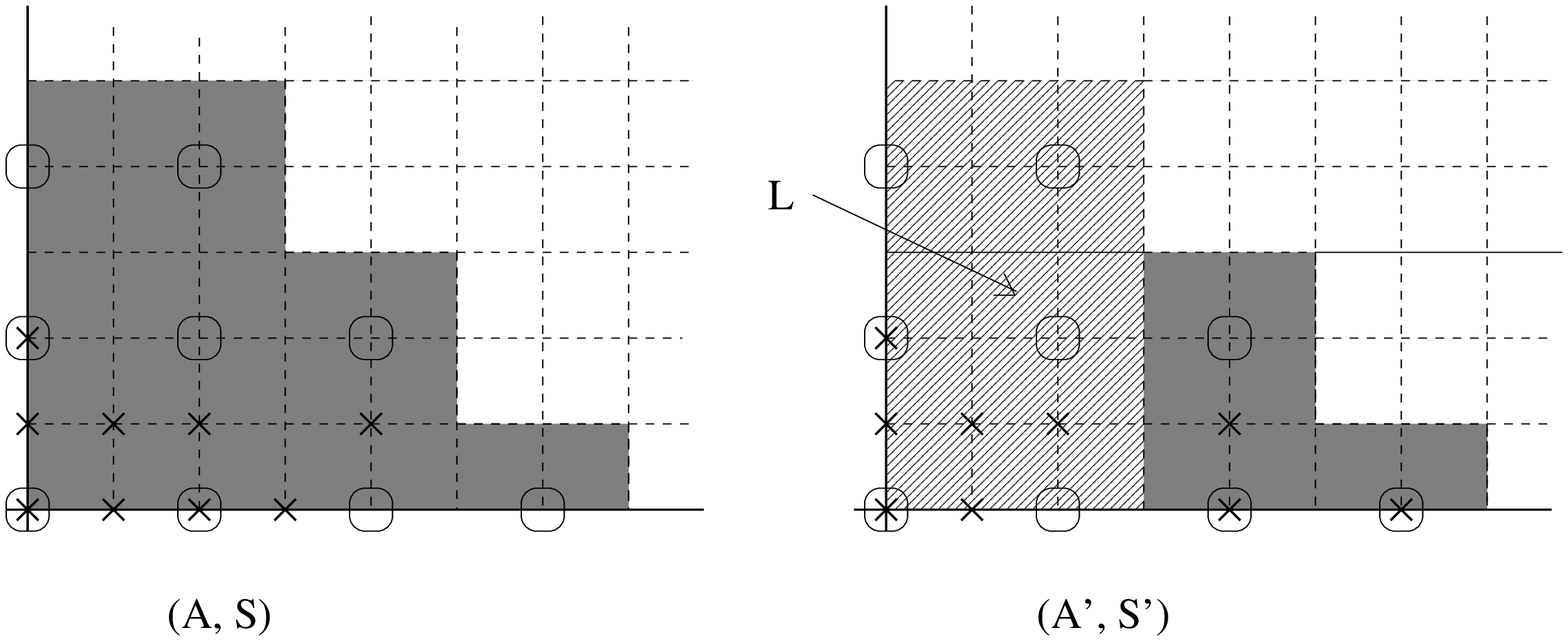}
\end{center}
\caption{}
\label{nonreg2}
\end{figure} 
\end{example}

% \begin{figure}[h]
% \begin{center}
%         \setlength\epsfxsize{7cm}
%         \leavevmode
%         \epsfbox{conreg2c.eps}
% \end{center}
% \caption{}
% \label{nonregc2}
% \end{figure} 
% \end{example}

\begin{example}\rm \
Let us return to the set $A$ appearing in the last example in \ref{half}
(Fig. ~\ref{nonreg2}),
and look for all $S$'s which make $(A, S)$ regular. As in Example \ref{ten4}, there are two possible cases.
One of them is precisely the one treated in Example \ref{half}, while the other one is shown in 
Fig. ~\ref{reg7}. The last one is regular. To see this, one first 
uses Theorem \ref{replace} to remove the last point of $A$ situated on $OY$.
The resulting scheme $(A', S')$ has $S_{x}(A)$ rectangular, hence we can use
Theorem \ref{rectangular} (or, more precisely, the version obtained by interchanging
$x$ and $y$). The subsequent univariate schemes
are regular (by \ref{Polya} again). 
\begin{figure}[h]
\begin{center}
         \setlength\epsfxsize{12cm}
         \leavevmode
         \epsfbox{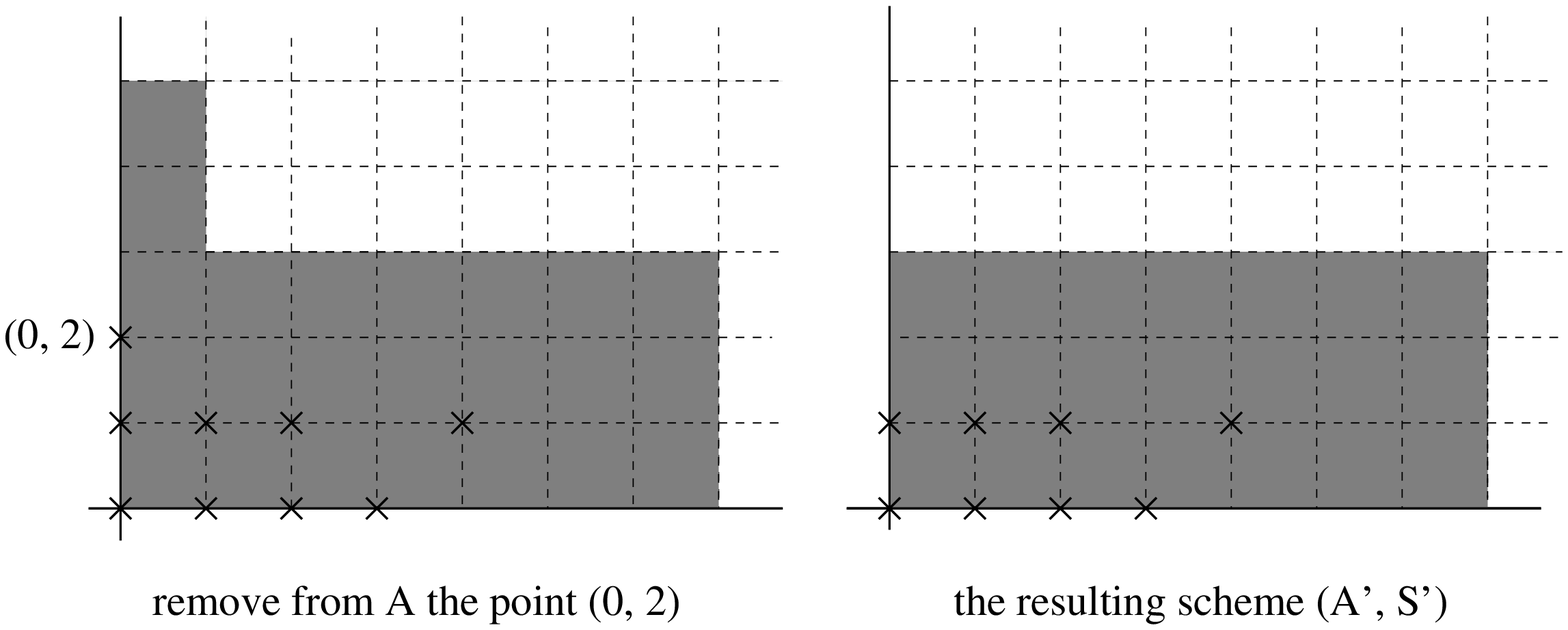}
 \end{center}
 \caption{}
 \label{reg7}
 \end{figure} 
\end{example}

\begin{example}\rm \
 A similar example is obtained by considering $A$ as in Fig. ~\ref{reg5}.
 As before, $S$ must be obtained from $S_0$ in the picture by adding
 one copy of $R(1, 1)$. There are two possible ways to do that, but 
 only the one shown in the picture produces a regular scheme. However, we do not
 know how to use the general results of the previous sections to prove the non-regularity
 of the other scheme.
 \begin{figure}[h]
 \begin{center}
         \setlength\epsfxsize{8cm}
         \leavevmode
         \epsfbox{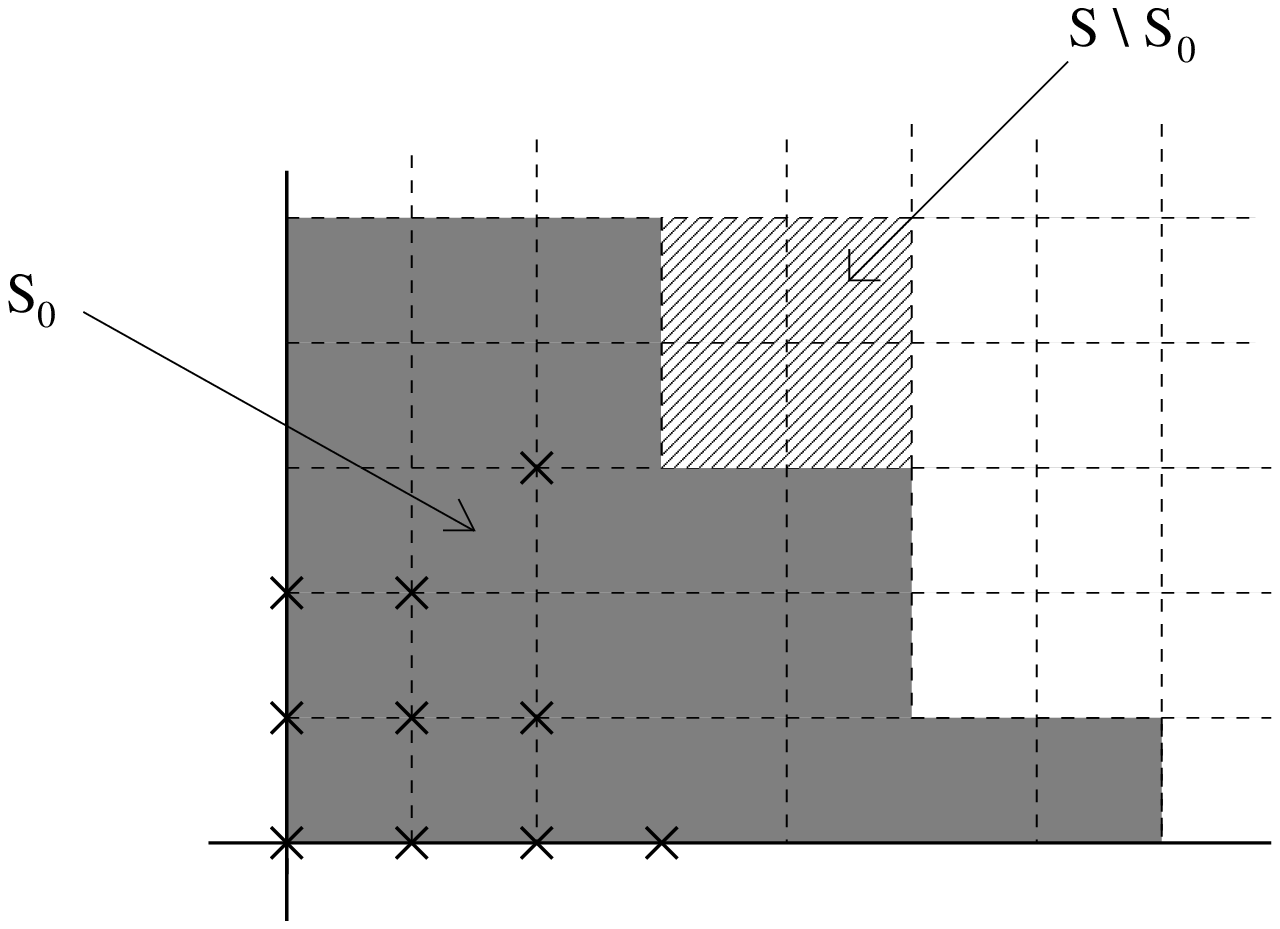}
 \end{center}
 \caption{}
 \label{reg5}
 \end{figure} 
 \end{example}

\begin{example}\rm \
Consider now $A$ as in Fig. ~\ref{reg3}.
As before, Theorem \ref{axes} tells us what $S$ must be around the axes.
We then have to fit three more (disjoint) copies of $R(1, 1)$
inside $R(5, 5)$ to get the lower set $S$. This time, this is possible in only one way
(as in the picture). One sees that the inequalities (\ref{str-Polya}) are satisfied, and one can actually show that the scheme is regular.
 \begin{figure}[h]
 \begin{center}
         \setlength\epsfxsize{6cm}
         \leavevmode
         \epsfbox{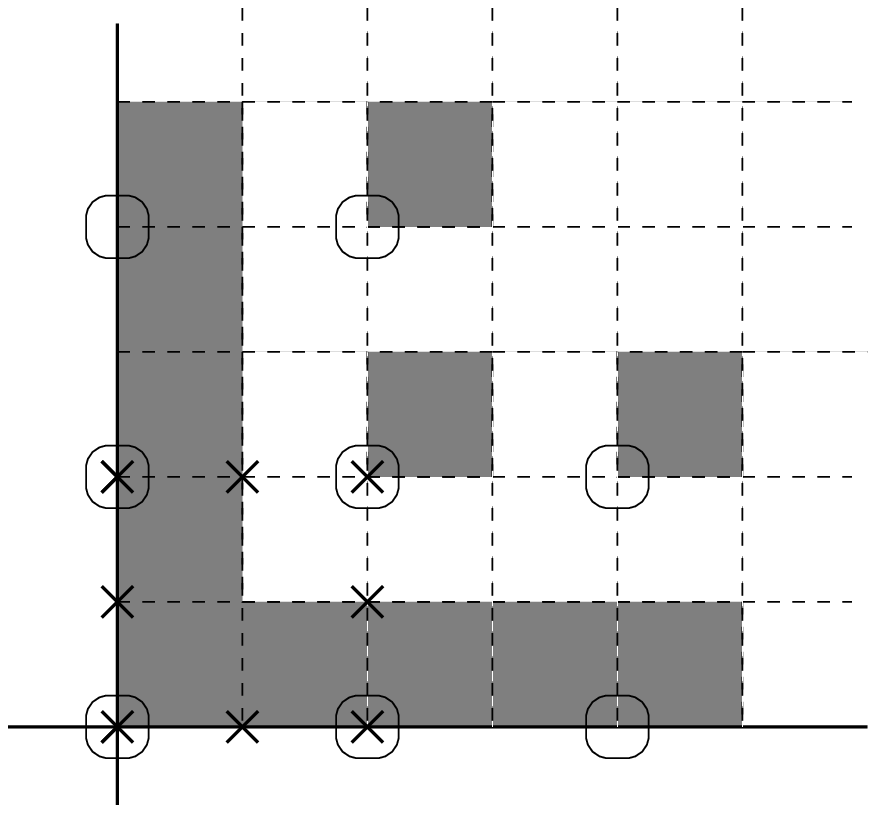}
 \end{center}
 \caption{}
 \label{reg3}
 \end{figure} 
\end{example}

\begin{figure}[h]
\begin{center}
        \setlength\epsfxsize{11cm}
        \leavevmode
        \epsfbox{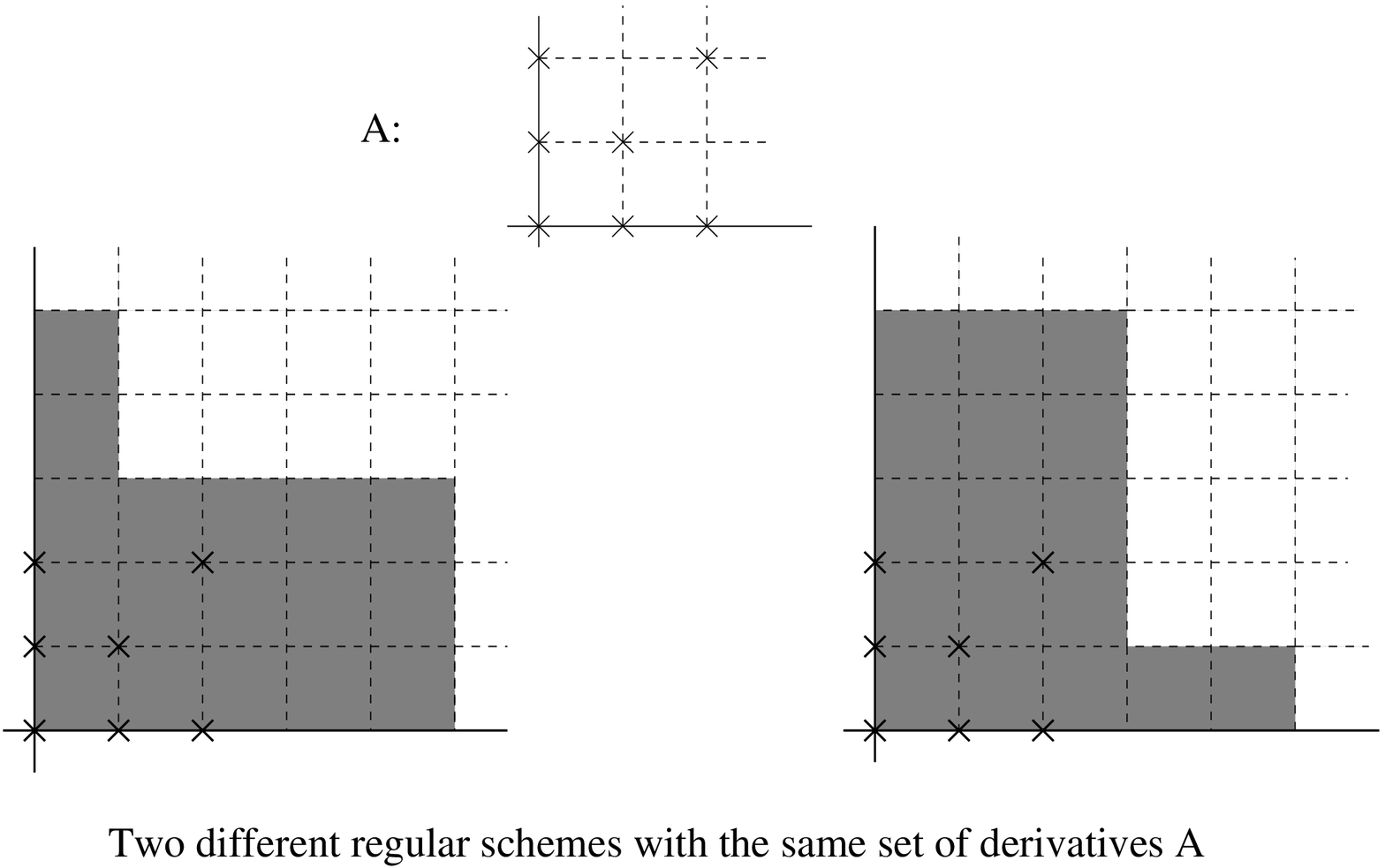}
\end{center}
\caption{}
\label{last2}
\end{figure}

\begin{example}\rm \
\label{two-of-them}
In all the previous examples, given $A$, there was at most one $S$ making 
$(A, S)$ regular. Here is an example where two such $S$'s can be chosen.
Consider $A$ as in Fig. ~\ref{last2}. As in the previous examples, $S$ is contained
in $R(5, 5)$ and must be obtained by adding one copy of $R(1, 1)$ to the blow
up $(T_{2})^{1, 1}$ of the triangle $T_{2}$. This is possible in two ways,
with the resulting $S$'s:  $S_{y}(2, 1, 1)^{1, 1}$ and $S_{y}(2, 2, 1)^{1, 1}$
(see also Fig. ~\ref{six}). That both resulting schemes are regular follows
again by removing one point and then using Theorem \ref{rectangular}. 

\begin{figure}[h]
\begin{center}
         \setlength\epsfxsize{11cm}
        \leavevmode
        \epsfbox{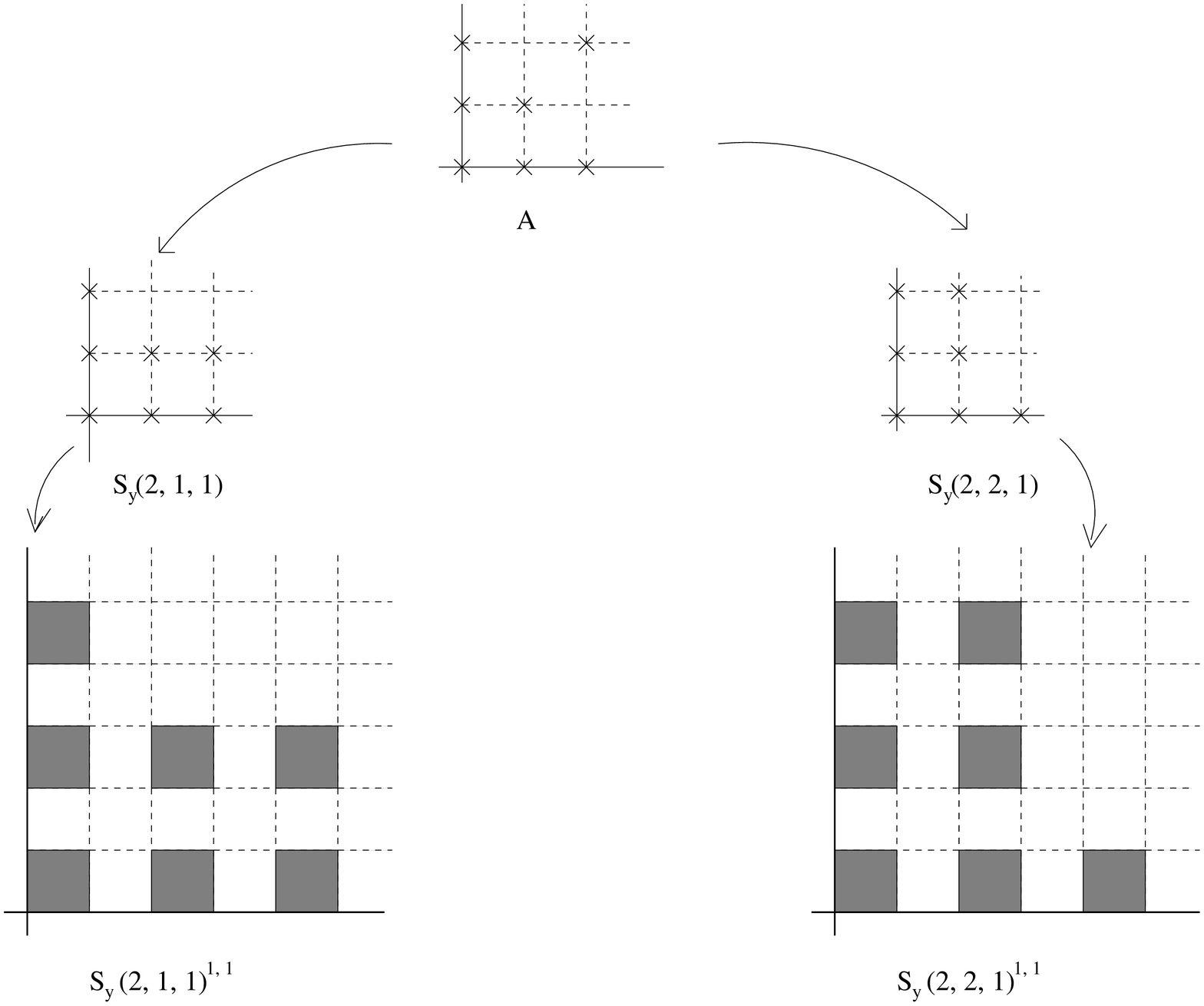}
\end{center}
\caption{First move A backwards to lower sets, and then blow up.}
\label{six}
\end{figure}
\end{example}

It is interesting to discuss this example in the light of Conjecture \ref{conj4}.
Note that, in all the examples we have considered so far, moving $A$ backwards to a lower set was possible in only one way,
hence the conjecture suggest the uniqueness of $S$ (proven by us in each case separately). 
In the example under discussion, the process of ``moving $A$ backwards'' to a lower set is not unique; there
are clearly (only) two ways of doing so, as shown in Fig. ~\ref{six}: the extremal element of $A$ can be moved one step down, or one step
to the left. And this is how we actually constructed this example. Similarly, one can find examples where 
the number of choices for $S$ equals a given number.

% \end{num}

% \begin{example}\rm \
% Using the results of the previous sections one can describe all the possible
% schemes $(Z, A, S)$ which are regular 
% with respect to $(p, q)$-rectangular sets of nodes, and which involve at most
% five derivatives (i.e. $|A|\leq 5$). Details will be given elsewhere.
% \end{example}

 % -----------------------------------------------------------------------
 \bibliographystyle{amsplain}
 \def\lllll{}
 
 \end{document}